\documentclass{amsproc}
\usepackage{amssymb, amsmath}

\newtheorem{theorem}{Theorem}[section]
\newtheorem{lemma}[theorem]{Lemma}
\newtheorem{corollary}[theorem]{Corollary}
\newtheorem{proposition}[theorem]{Proposition}

\theoremstyle{definition}
\newtheorem{definition}[theorem]{Definition}
\newtheorem{example}[theorem]{Example}

\theoremstyle{remark}
\newtheorem{remark}[theorem]{Remark}
\newtheorem*{acknowledgments}{Acknowledgments}

\numberwithin{equation}{section}

\begin{document}
\title[Construction of Parseval wavelets from redundant filter systems]
{Construction of Parseval wavelets from redundant filter systems}

\author{L. W. Baggett}
\address{Department of Mathematics, Campus Box 395, University of Colorado, Boulder, CO, 80309-0395}
\email{baggett@euclid.colorado.edu}
\author{P. E. T. Jorgensen}
\address{Department of Mathematics, University of Iowa, 14 MacLean Hall, 
Iowa City, IA, 52242-1419}
\email{jorgen@math.uiowa.edu}
\thanks{The first two named authors were supported by a US-NSF 
Focused Research Group (FRG) grant.}
\author{K. D. Merrill}
\address{Department of Mathematics, Colorado College, Colorado Springs, CO 80903-3294}
\email{kmerrill@coloradocollege.edu}
\author{J. A. Packer}
\address{Department of Mathematics, Campus Box 395, University of Colorado, Boulder, CO, 80309-0395}
\email{packer@euclid.colorado.edu}
\date{June 9, 2003}

\begin{abstract}
We consider wavelets in $L^2(\mathbb{R}^d)$ which have generalized multiresolutions.
This means that the initial resolution subspace $V_0$ in $L^2(\mathbb{R}^d)$ is not singly generated.
As a result, the representation of the integer lattice $\mathbb{Z}^d$
restricted to $V_0$ has a nontrivial multiplicity function.
We show how the corresponding analysis and synthesis for these wavelets
can be understood in terms of unitary-matrix-valued functions on a torus acting
on a certain vector bundle. Specifically, we show how the wavelet functions on
$\mathbb{R}^d$ can be constructed directly from the generalized wavelet filters.
\end{abstract}
\maketitle


\section{Introduction}

The theory of wavelets is concerned with the Hilbert space $L^{2}\left(
\mathbb{R}^{d}\right)  $. The problem is to find ``good'' orthonormal bases
(ONB), where ``good'' makes reference to several conflicting requirements:

\begin{enumerate}
\item \label{Req1}These bases must be constructed from a small number of model
functions, called wavelets, and two discrete operations, translation and
scaling. In this paper, we are concerned with translation by the standard
integer lattice $\mathbb{Z}^{d}$, and scaling by some prescribed integral
matrix $A$ which is assumed expansive.

\item \label{Req2}In passing from function to expansion coefficients,
referring to a wavelet basis, and back again (this is called
analysis/synthesis), the steps must be algorithmic, ideally avoiding direct
reference to integration over $\mathbb{R}^{d}$.

\item \label{Req3}The wavelet functions should have compact support, and
should have some prescribed number of derivatives.
\end{enumerate}

The algorithms that have been popular since the mid 1980's are based on what
is called multiresolution analysis (MRA). This was pioneered by Daubechies
\cite{D}, Mallat \cite{Ma} and Meyer \cite{Me,Mey93b}, and the idea, while
simple, has been extremely powerful. The idea itself is much like that of the
Gram-Schmidt algorithm from Hilbert space, in that it is based on a scale of
closed subspaces, resolutions $V_{n}$, and relative orthogonal complements,
detail subspaces $W_{n}$. The scale of subspaces $V_{n}$ play the role of
martingales from probability theory.

Daubechies's book \cite{D} stresses how the requirements (\ref{Req1}%
)--(\ref{Req3}) can be met with the MRA approach, and all starting with a
fixed cleverly chosen function $\varphi$ in the subspace $V_{0}$ from the
resolution. The function $\varphi$, the father function, is the solution in
$L^{2}$ to a scaling, or refinement, equation; a solution which results from a
cascade approximation. The wavelet functions, mother functions, can then be
constructed from the subspace $W_{0}$ which is the relative ortho-complement
of $V_{0}$ in $V_{1}$. One drawback of this approach is that if $N=\left\vert
\det A\right\vert $, then $N-1$ wavelet generators from $W_{0}$ are needed. In
the dyadic case, $d=1$, $N=2$, that makes one function, but in general, $N$
can be large. Now the spaces $V_{0}$ and $W_{0}$ are invariant under
translation by $\mathbb{Z}^{d}$, and there is a corresponding pair of
multiplicity functions which dictate a minimal choice of generators for
$V_{0}$ and $W_{0}$.

In fact for the general case of $L^{2}\left(  \mathbb{R}^{d}\right)  $ and a
fixed scaling matrix $A$, it is possibly to get $W_{0}$ singly generated,
i.e., to find a single generator $\psi$. In some cases $\psi$ may be taken to
be the inverse Fourier transform of a subset $E$ of $\mathbb{R}^{d}$. Such
subsets $E$ are called wavelet sets, see \cite{DaLa98,DLS97,DLS98} and \cite{BJMP,BMM,BM}%
. But there are other choices of sets of generating functions $\psi$, with the
number of generators between $1$ and $N-1$. What emerges is that these
constructions force frequency localization, and the compact support in the
$x$-variable is typically lost, i.e., we must relax requirement (\ref{Req3}).
As it turns out, the kind of frequency localization we obtain is well suited
for effective sampling algorithms. A second issue enters: The wavelet
algorithm may lead to bases which only satisfy a certain Parseval property
(also called ``normalized tight frame''). While we still have the resolution
structure $V_{n}$, the number of generators in $V_{0}$ may increase, but they
are not directly part of the wavelet basis. This setup is referred to as a
generalized multiresolution analysis (GMRA); see \cite{BJMP,BMM,BM}.

The GMRA theory was in fact introduced (in an operator-theoretic/operator-%
algebraic framework) in a pioneering paper by Baggett, Carey, Moran, and
Ohring \cite{BCMO95}.
By now there is a rich journal literature which reflects wavelet
constructions with some degree of multiplicity; see, e.g.,
\cite{ACM04,BL,DHRS03,Han97,HLPS,KiLi01,Li01,PSWX01,Pap00,RoSh97,Sel01,Web02,Web04}.
In addition our multiplicity analysis has
recently been used in papers by Dutkay and Jorgensen
\cite{DuJo03,DuJo04a,DuJo04b,DuJo04c}
in
the study of nonlinear dynamical models.

       A recent and interesting paper of Papadakis \cite{Pa04} offers a global and natural approach to generalized multiresolution analyses (GMRAs)  based on a geometric frame construction which
has the advantage of including all GMRA constructions in $L^2(\mathbb{R})$. Some of the differences between the approach in \cite{Pa04} and the present one lie in our use here of operator algebraic tools
deriving from the Cuntz algebras \cite{Cun} in operator theory. Further, our use of vector bundles in Section 4 offers a rather explicit and concrete representation of the matrix functions which in turn
describe the filters for GMRA normalized tight frame wavelets, alias, for Parseval wavelet frames.

       Our paper and other recent papers on wavelet frames may be said to generalize a celebrated theorem first proved by Lawton in this
journal in the paper \cite{L}.  Lawton's pioneering result states that a trigonometric low-pass filter which satisfies a certain
conjugate-mirror filter condition must give rise to a wavelet frame. However, the resulting wavelet frame may not be associated with a
classical MRA. We generalize Lawton's theorem in several directions to a much broader class of subband wavelet filters.
Our matrix-valued low-pass filters are still associated with GMRA wavelets. These
in turn are the most general types of multiresolution structures. 

    Finally, we emphasize that the use of wavelet filters derives historically from signal processing in communications engineering, see
e.g., \cite{Jor03}.

In the standard case of MRAs, it is well known how the subband filters from
signal processing allow us to construct the wavelet functions by an elegant
algorithm. The function $\varphi$ is a solution to a certain refinement
equation. The subband filters may be thought of as functions on a torus
$\mathbb{T}$, frequency response functions. But in the case of multiplicity
and multiple generators, the corresponding functions on $\mathbb{T}$ are
matrix-valued, and the refinement equation is a matrix equation.

In this paper we show that the generalized setup admits solutions in
$L^{2}\left(  \mathbb{R}^{d}\right)  $ starting with this matrix/vector
version of the refinement equation. Starting with a matrix system of subband
filter functions on a torus, we show that our corresponding wavelet solutions
are in $L^{2}\left(  \mathbb{R}^{d}\right)  $, and that they will be Parseval
frames for the Hilbert space $L^{2}\left(  \mathbb{R}^{d}\right)  $.

While this is a weakening of the stricter ONB in requirement (\ref{Req2}), the
Parseval frame property still allows the same recursive analysis/synthesis
algorithms 
popular in the MRA case.
The now classical method of Mallat and Meyer for constructing an orthonormal wavelet
in $L^2({\mathbb R})$ (relative to translation by integers and dilation by 2
 \cite{Ma,Me})
proceeds as follows:  Let $h$ be a periodic function on ${\mathbb R}$ that satisfies the
``conjugate mirror filter condition'' (the so-called Smith-Barnwell condition)
\begin{equation}
\label{classicalorthh} 
|h(x)|^2+|h(x+\frac12)|^2 =2.
\end{equation}
The function $h$ is the conjugate mirror filter, referred to above, and in our context is   
called a \textit {low-pass filter.}
Consider the infinite product 
$$P(x)=\prod_{j=1}^\infty \frac1{\sqrt2} h(\frac x{2^j}),$$
and suppose that there exists a nonzero $L^2$ function $\phi$ whose Fourier transform $\widehat\phi$
coincides with $P.$
Under not too strenuous assumptions on $h,$
this does in fact hold.
For instance, if we take $h$ to be smooth, 
and satisfying Cohen's orthogonality conditions and the \textit{low-pass condition}
$|h(0)|=\sqrt2$,
the set of integer translates of $\phi$ turn out to be orthonormal functions in $L^2$, 
and in fact $\phi$ is a scaling function for a multiresolution analysis $\{V_j\}.$
\par
Given such a low-pass filter function $h,$ there exists an associated periodic
function $g,$ which also satisfies the filter equation \ref{classicalorthh}, and such that
$h$ and $g$ satisfy the following orthogonality condition:
\begin{equation}
\label{classicalorthgh}
h(x)\overline {g(x)} + h(x+\frac12)\overline{g(x+\frac12)} = 0.
\end{equation}
Any such function $g$, called a \textit{high-pass filter}, can be obtained from the low-pass fil-\linebreak ter
function via the standard technique of constructing a unitary matrix
whose first row is given by $((h(x)/\sqrt2),(h(x+1/2)/\sqrt2).$
Finally, the function
$\psi,$ defined  by
$$\widehat\psi(x) =\frac1{\sqrt 2} g(\frac x2)\widehat\phi(\frac x2)$$ 
is an orthonormal wavelet.  That is, the collection
$\{\psi_{j,n}\}\equiv \{\sqrt2^j\psi(2^jx-n)\},$ for $j$ and $n$ in ${\mathbb Z},$ forms an orthonormal
basis for $L^2({\mathbb R}).$
\par
A famous example of A. Cohen (\cite{Ch}) shows that eliminating the non-vanishing condition 
can cause the Mallat-Meyer method to go wrong in an interesting but not disastrous way.  
Cohen exhibited a low-pass filter function $h,$ for which the infinite product
$P$ exists, is the Fourier transform of a nonzero $L^2$ function $\phi,$ but for which the
integer translates of $\phi$ are not orthonormal.  Further, the
translates and dilates $\{\psi_{j,n}\}$ of the corresponding function $\psi,$
defined just as in the classical method, is  not an orthonormal basis.
So, for this choice of filter $h,$ the Mallat-Meyer procedure
fails to produce a scaling function, and the resulting function $\psi$ is not an
orthonormal wavelet.  Nevertheless, its translates
and dilates do form what's called a {\sl Parseval frame.\/}
By definition, this means that
for each $f\in L^2({\mathbb R})$ we have 
$$\|f\|^2 = \sum_j\sum_n |\langle f\mid \psi_{j,n}\rangle|^2.$$
\par
In \cite{L} and \cite{BJ}, the Mallat-Meyer phenomenon was generalized to incorporate the Cohen
example in the following way.
Suppose $h$ is a low-pass filter for dilation by a positive integer $N,$ i.e., satisfies the filter equation  
$$\sum_{l=0}^{N-1} |h(x+\frac lN)|^2 = N$$
and the low-pass condition $|h(0)|=\sqrt N,$
and suppose $g_1,\ldots,g_{N-1}$ are corresponding high-pass filters, i.e., 
are periodic functions for which
the $N$ functions $h,g_1,\ldots,g_{N-1}$ satisfy the following orthogonality conditions:
$$\sum_{l=0}^{N-1} h(x+\frac lN) \overline{g_i}(x+\frac lN) = 0$$
for all $1\leq i \leq N-1,$ and
$$\sum_{l=0}^{N-1} g_i(x+\frac lN)\overline {g_j(x+\frac lN)}
=\begin{cases} N & i=j\\
0& i\neq j
\end{cases},$$
for all $i,j$ between 1 and $N-1.$
Again, these $N-1$ high-pass filters can be constructed from the low-pass
filter $h$ by the matrix completion technique.
\par
Let $P$ be the infinite product
$$P=\prod_{j=1}^\infty \frac1{\sqrt N} h(\frac  x{N^j}).$$
Then, if the functions $h,g_1,\ldots,g_{N-1}$ are Lipschitz continuous,
it is shown in \cite{BJ}
that there exists a nonzero $L^2$ function $\phi$
whose Fourier transform coincides with $P,$ 
and the $N-1$ functions
$\{\psi_k\},$ defined by
$$\widehat{\psi_k}(x ) =\frac1{\sqrt N} g_k(\frac x  N)\widehat\phi(\frac x  N),$$
form a Parseval frame multiwavelet.  That is, the collection
$$\{\psi_{j,n,k} \equiv \sqrt N^j\psi_k(N^jx-n)\},$$
 $j,n\in{\mathbb Z}$ and $1\leq k \leq N-1,$
forms a Parseval frame for $L^2({\mathbb R}).$
\par
Just as in the Cohen phenomenon in the Mallat-Meyer constructions, the integer translates
of the function $\phi$ may or may not be orthonormal,
even though their closed linear span $V_0$ does form the core subspace of a (generalized) multiresolution analysis.
Also, the wavelets $\psi_k$ may or may not have orthonormal translates and dilates, and may or may not belong to the subspace $V_1\ominus V_0$ of this associated GMRA. 
\par
The arguments in the proof of this result introduce some
ideas from operator theory, specifically in the form
of the Ruelle operator and partial isometries that satisfy relations similar to the Cuntz relations \cite{Cun}.  In the case of dilation by 2, W. Lawton \cite{L} used a cascade algorithm in the time domain to independently derive the same result
(in the special case of trigonometric polynomials, i.e., the case corresponding to a single compactly supported scaling function).
\par
While the constructions of \cite{BJ}and \cite{L} start with filters associated with a classical multiresolution analysis, they build a frame wavelet that may be obtained from only a generalized multiresolution analysis.   
The purpose of the present paper is to extend and clarify this kind of result to filters defined in higher dimensional and non MRA contexts.  That is, we suppose that $A$ is an expansive,  integral, $d\times d$ matrix, 
and we investigate frame wavelets constructed
from filter systems associated to generalized multiresolution analyses in $L^2({\mathbb R}^d)$ relative
to dilation by $A$ and translation by lattice points in ${\mathbb Z}^d.$
\par
The theory of generalized conjugate mirror filters relative to
a generalized multiresolution analysis was first developed in \cite{Cou} and \cite{BCM}.
In Section 2 we
briefly review that subject, and in particular we recall the
analogs to Equations \ref{classicalorthh} and \ref{classicalorthgh}, i.e., generalized filter and
orthogonality equations.
These generalized equations are considerably different from the
original ones, because the right-hand side now involves an integer-valued multiplicity function on
the $d$-dimensional torus ${\mathbb T}^d\equiv [-\frac 12,\frac 12)^d.$
In general, the low-pass filter becomes a square matrix $H=[h_{i,j}]$ of periodic functions, and the
high-pass filter becomes a not necessarily square matrix $G=[g_{k,j}]$ of periodic functions.
Because the filters in this generalized context are matrices of periodic
functions, we refer to them as redundant systems of filters.
\par
By generalizing the arguments in \cite{BJ}, we are able here to build wavelets from 
generalized filters more simply and with fewer restrictions than in \cite{BCM}.
In Section 3, under smoothness and low-pass conditions analogous to those in \cite{BJ},
we prove that the infinite matrix product
$$P = \prod_{j=1}^\infty \frac1{\sqrt{|\det(A)|}} H({A^t}^{-j}(x))$$
converges pointwise, and the first column $\Phi$ of $P$
is the Fourier transform of a set $\phi_1,\phi_2,\ldots$ of $L^2$ functions that together generate a generalized multiresolution analysis.
We then prove that
the functions $\{\psi_k\}$ given by
$$\widehat{\psi_k}(x) = \frac1{\sqrt{|\det(A)|}} \sum_j g_{k,j}({A^t}^{-1}(x)) \widehat{\phi_j}({A^t}^{-1}(x))$$
form a Parseval wavelet frame for $L^2({\mathbb R}^d).$
\par In engineering, the various classes of multiresolution analyses are
motivated by filters from signal and image processing. While these are
practical concerns, our present aim is theoretical. In most of the
engineering MRA-constructions, the multiplicity function is equal to $1$ on
the entire fundamental domain $[-1/2,1/2)$, so this only produces
classical MRA quadrature-mirror filters.  One must choose the multiplicity
function to be mixed, i.e., equal to the characteristic function of
delicately selected configurations of subsets of $[-1/2,1/2)$, to get
genuinely generalized filters, i.e., GMRA filters. If the multiplicity function is not constant, the filters can never all be continuous in the frequency domain, which can be a problem for engineers.  The reason that all the filters for a GMRA cannot be continuous is that they must satisfy the filter equations \ref{orthh} and \ref{orthg}, both of which involve non-continuous characteristic functions for non-constant multiplicity functions.
This particular feature of our theory becomes evident in 
our example in Section 3, Example \ref{ExJourne}. Yet,
this example is surprising in two ways: 
Firstly, it is related to a non-classical MRA-wavelet, the Journ\'e wavelet. 
Secondly, we construct two low-pass filters which have discontinuities in
the frequency domain, but the resulting scaling functions are $C^\infty$ in
the frequency domain! This prompts the following question: can the scaling functions and wavelets associated to GMRA's 
be constructed to have nice properties in both the time and the frequency domains? From the
viewpoint of engineering , this is an important question, and it is 
addressed in a sequel paper \cite{BJMP05} in much more detail. We
have established this in the affirmative in  \cite{BJMP05}, where we construct a frame wavelet that is $C^r$ in the time domain and $C^{\infty}$ in the frequency domain.  More importantly, it would be
interesting to know whether or not scaling functions and frame wavelets that have nice properties in both the time and frequency domains can be
constructed for every choice of a multiplicity function. 
\par
In \cite{BJ}
the set of functions $h,g_1,\ldots,g_{N-1}$ is called
an \textit{M-system.} It is shown in \cite{BJ} that there is
a group that acts freely and transitively on the set of $M$-systems, thus suggesting a natural
structure on these systems and therefore on the corresponding frame wavelet systems.
In Section 4 of this paper, we describe an analogous action on the generalized filter systems we have introduced.  This time,
it is a group bundle that acts freely and transitively.


\section{Generalized filters}
We collect here the relevant definitions concerning wavelets and multiresolution analyses
in $L^2({\mathbb R}^d),$ relative to translation by lattice points and dilation
by an expansive integer matrix $A;$
i.e., a matrix each of whose eigenvalues has modulus greater than 1.
\par
Recall that a {\sl frame\/} in $L^2({\mathbb R}^d)$ is a sequence $\{f_n\}$ for which there exist positive numbers $a$ and $b$ such that
$$a\|f\|_2^2 \leq \sum_n |\langle f\mid f_n\rangle|^2 \leq b\|f\|_2^2$$
for every $f \in L^2({\mathbb R}^d).$
It is called a {\sl Parseval frame\/} or a {\sl normalized tight frame\/} if
$a=b=1$ in the inequalities above.  That is, $\{f_n\}$ is a Parseval frame if
$\|f\|^2_2 = \sum_n |\langle f\mid f_n\rangle|^2$
for every $f\in L^2({\mathbb R}^d).$
\par
For each $z\in {\mathbb Z}^d,$ we write $\tau_z$ for the unitary {\sl translation operator\/} on $L^2({\mathbb R}^d)$ given by
$[\tau_z(f)](t) = f(t+z).$
Fix an expansive, $d\times d,$ integer matrix $A$ and let $B=A^t$ and $N=|\det(A)|.$
We write $\delta$ for the unitary {\sl dilation operator\/} given by
$[\delta(f)](t) = \sqrt{ N}f(A(t)).$  For each element $\omega\in {\mathbb T}^d,$ there exist exactly $N$
distinct points $\zeta\in {\mathbb T}^d$ such that $\alpha(\zeta)=\omega,$ where
$\alpha$ denotes the endomorphism on ${\mathbb T}^d$ induced by the action of $B$ on ${\mathbb R}^d.$
We denote these points $\zeta$ in a Borel manner as $\omega_0,\omega_1,\ldots,\omega_{N-1}.$
Note also that, because $A$ is an expansive matrix,
the endomorphism $\alpha$ is ergodic.  

\begin{definition}
\label{GMRAdef}
A \textit{multiresolution analysis} (MRA) in $L^2({\mathbb R}^d),$
relative to the group $\{\tau_z\}$ of translation operators and the dilation operator $\delta,$
is a sequence $\{V_j\}_{-\infty}^\infty$ of closed subspaces of $L^2({\mathbb R}^d)$ for which:
\begin{enumerate}
\item\hskip0.5em $V_j\subseteq V_{j+1}.$
\item\hskip0.5em $V_{j+1} = \delta(V_j).$
\item\hskip0.5em $\cup V_j$ is dense in $L^2({\mathbb R}^d),$ and
$\cap V_j = \{0\}.$
\item  There exists an element $\phi\in V_0$ whose translates $\{\tau_z(\phi)\}$ form an orthonormal basis for $V_0.$
Such an element $\phi$ is called a {\sl scaling function\/} for the multiresolution analysis.
\end{enumerate}
\par
A \textit {generalized multiresolution analysis} (GMRA) is a sequence $\{V_j\}$ of closed
subspaces of $L^2({\mathbb R}^d)$ that satisfies conditions (1), (2), and (3) above,
but satisfies the weaker fourth condition
\begin{enumerate}
\item[$4'.$]\hskip0.5em $V_0$ is invariant under all the translation operators $\tau_z.$
\end{enumerate}
In both cases, the subspace $V_0$ is called the {\sl core subspace\/} of the GMRA.
\end{definition}
The theory of GMRA's is developed in \cite{BCMO95}, \cite{BMM}, and \cite{BM}
(see also \cite{Pap00}).  In particular, it is shown there that by Stone's Theorem for unitary representations of the
group $\mathbb Z^d$, or more generally by
spectral multiplicity theory for a set of commuting unitary operators, the invariance under translation of the core subspace $V_0$ implies the existence of a \textit{multiplicity function}
$m:\mathbb T^d\mapsto\{0,1,2,\cdots\infty\}.$  The multiplicity function $m$ counts the number of times each exponential function is represented as a subrepresentation of 
translation by $\mathbb Z^d$ on $V_0$.  We will assume in this paper that $m$ is bounded, with maximum value $c$.  We let $S_i=\{\omega\in\mathbb T^d:m(\omega)\geq i\},$ and recall from \cite{BM} that there exist \textit{generalized scaling functions}
$\{\phi_i\}_{1\leq i\leq c}$ such that the collection $\{\tau_z(\phi_i)\}$
for $z\in{\mathbb Z}^d$  and $1\leq i\leq c$ forms a Parseval frame for $V_0,$ and such that
$\sum_{z\in\mathbb Z^d}|\widehat{\phi}_i(\omega + z)|^2=\chi_{S_i}(\omega)$.  Note that these results translate to 
the classical conditions in an MRA, which is the special case of a GMRA in which the multiplicity function is identically 1.

If $V_0$ is the core subspace of a generalized multiresolution analysis $\{V_j\},$
then the subspace $W_0=V_1\cap V_0^\perp$ of a GMRA $\{V_j\}$ also is
invariant under all the translation operators $\tau_z.$  Hence, again by Stone's theorem,
there exists a \textit{complementary multiplicity function} $\widetilde{m}:\mathbb T^d\mapsto\{0,1,2,\cdots\infty\}$ that 
characterizes the representation of $\mathbb Z^d$ on $W_0$.  
As a direct result of the fact that $V_1=V_0\oplus W_0$, the multiplicity functions $m$ and $\widetilde{m}$ must satisfy
the following \textit{consistency equation} (see \cite{BMM}):
\begin{equation}
\label{consistency}
m(\omega) + \widetilde m(\omega) = \sum_{l=0}^{N-1} m(\omega_l).
\end{equation}
By the consistency equation, the assumption that $m$ is bounded implies that $\widetilde{m}$ is bounded as well.  We write
$\widetilde{c}$ for the maximum value of $\widetilde{m}$, and $\widetilde{S}_k$ for $\{\omega\in\mathbb T^d: \widetilde{m}(\omega)\geq k\}.$

Generalized multiresolution analyses are useful because of their relationship to wavelets.  In particular, in \cite{BMM} it is shown that every 
orthonormal multiwavelet is associated to a GMRA with $\widetilde{m}=$the number of wavelets. We recall the relevant definitions here:   
\begin{definition}
\label{waveletdef}
An \textit{orthonormal wavelet,} or more generally
 an \textit{orthonormal multiwavelet,} in $L^2({\mathbb R}^d),$
relative to the group $\{\tau_z\}$ of translation operators and the dilation operator $\delta,$
is a finite set $\psi_1,\ldots,\psi_m$ of functions in $L^2({\mathbb R}^d)$
for which the collection $\{\psi_{j,z,k}\}\equiv \{\delta^j(\tau_z(\psi_k))\}$
forms an orthonormal basis for $L^2({\mathbb R}^d).$
A set of functions $\psi_1,\ldots,\psi_m$ forms a \textit{frame multiwavelet} if the collection $\{\psi_{j,z,k}\}$ forms a frame for $L^2({\mathbb R}^d),$
and the set is called a \textit{Parseval frame multiwavelet} if the collection $\{\psi_{j,z,k}\}$
forms a Parseval frame for $L^2({\mathbb R}^d).$
\end{definition}

We now are ready to develop the definition of generalized filters.  The classical filter equation (Equation \ref{classicalorthh}), which is 
the basis for the Mallat-Meyer construction of scaling functions and wavelets,
 emerges naturally from
the study of classical multiresolution analyses.  For, if $\phi$ is a function, for which 
the closure of the linear span of the translates of $\phi$ is invariant
under the dilation
$f(x) \to f(x/2),$ 
e.g., if $\phi$ is a scaling function for an MRA,
then it follows from elementary Fourier analysis that $\widehat\phi$ satisfies the following classical \textit{refinement equation:}
\begin{equation}
\label{classrefeq}
\sqrt2 \widehat\phi(2x) = h(x)\widehat\phi(x),
\end{equation}
where $h$ is a periodic function that satisfies the filter equation 1.1.  
As described in the introduction, the Mallat-Meyer construction \cite{Ma},\cite{Me} reverses this procedure, by beginning with a sufficiently well behaved 
function $h$ that satisfies filter equation 1.1, and producing a corresponding scaling function and multiresolution analysis by iterating equation \ref{classrefeq} to get an infinite product expression for 
$\widehat{\phi}$.

In an analogous way, the theory of generalized filters emerges naturally from the study of generalized multiresolution analyses (see \cite{Cou} and \cite{BCM}).  
Because $V_{-1} = \delta^{-1}(V_0)$ is contained in $V_0$, it follows
that, for each $i,$ there exists a sequence $\{h_{i,j}\}$ of functions on the torus ${\mathbb T}^d$ such that
\begin{equation}
\label{genrefeqn}
\sqrt N \widehat{\phi_i}(B(x)) = \sum_j h_{i,j}(x) \widehat{\phi_j}(x).
\end{equation}
As is shown in \cite{BCM}, these generalized filters satisfy the \textit{generalized filter equation}:
\begin{equation}
\label{orthh}
\sum_j \sum_{l=0}^{N-1} h_{i,j}(\omega_l)\overline{h_{i',j}(\omega_l)} =
\delta_{i,i'} N\chi_{S_i}(\omega),
\end{equation}
and also have $h_{i,j}$ supported on $S_j$.
In analogy with the classical case, we make the following definition.

\begin{definition}
\label{lowpassdef}
A \textit{generalized low-pass filter} relative to a GMRA with multiplicity functions $m$ and $\widetilde{m}$
is a matrix of functions $H=[h_{i,j}]$ on $\mathbb T^d$ (or equivalently $\mathbb Z^d$ periodic functions on
$\mathbb R^d$), with
$h_{i,j}$ supported on $S_j$ (or the periodization of $S_j$), that satisfy both the generalized filter equation (\ref{orthh})
and the generalized low-pass condition $h_{i,j}(0)=\delta_{i,1}\delta_{j,1}\sqrt N$
\end{definition}

Just as in the classical case, we can sometimes reverse the procedure of producing filters from wavelets.  
In particular, generalized low-pass filters give rise to generalized high-pass filters.  
The relevant result is contained in the following theorem, again from \cite{BCM}.
\begin{theorem}
Let $H=[h_{i,j}]$ be a generalized low-pass filter
relative to a GMRA with multiplicity functions $m$
 and  $\widetilde m.$
Then there exists a matrix $G=[g_{k,j}]$ of functions on ${\mathbb T}^d$ satisfying
\begin{equation}
\label{orthg}
\sum_j \sum_{l=0}^{N-1} g_{k,j}(\omega_l) \overline{g_{k',j}(\omega_l)} = 
\delta_{k,k'} N \chi_{\widetilde{S_k}}(\omega),
\end{equation}
and
\begin{equation}
\label{orthgh}
\sum_j \sum_{l=0}^{N-1} h_{i,j}(\omega_l)\overline{g_{k,j}(\omega_l)} = 0
\end{equation}
for all $i$ and $k.$
The matrix of functions $G=[g_{k,j}]$ is called a \textit{generalized high-pass filter}
\end{theorem}

Under conditions which allow the production of generalized scaling functions from generalized low-pass filters, these high-pass filters can
be used to build frame multiwavelets.  Some narrow conditions that allow this are described in \cite{BCM}.  In Section 3, we show that more
general conditions allow the same construction.  

There is one final property of a generalized low-pass and high-pass filter that we will need in Section 3,
which was not presented in \cite{BCM}.
\begin{theorem}
\label{thmorthocol}
Let $H$ and $G$ be a generalized low-pass and high-pass filter
relative to multiplicity and complementary multiplicity functions $m$ and $\widetilde m.$
Assume that the maximum value of $m$ is $c$ and the
maximum value of $\widetilde m$ is $\widetilde c.$
Then
\begin{equation}
\label{orthcol}
\sum_{i=1}^c  h_{i,j}(\omega_l)\overline{h_{i,j'}(\omega_{l'})}
+ \sum_{k=1}^{\widetilde c} g_{k,j}(\omega_l)\overline{g_{k,j'}(\omega_{l'})}
= \delta_{j,j'}\delta_{l,l'} N\chi_{S_j}(\omega_l).\end{equation}
\end{theorem}
\begin{proof}
For each $\omega\in{\mathbb T}^d,$ we define a matrix $K(\omega)$
having $c+\widetilde c$ rows and $C\times N$ columns as follows:
We index the $c\times N$ columns of $K(\omega)$ by pairs $(j,l),$ where $1\leq j \leq c$ and $0\leq l \leq N-1.$
Then the entry $k_{i,(j,l)}(\omega)$ is defined to be $(1/\sqrt N)h_{i,j}(\omega_l)$ if $1\leq i \leq c,$
and $k_{i,(j,l)}(\omega) =  (1/\sqrt N)g_{i-c,j}(\omega_l)$ if $c<i \leq \widetilde c.$
We see directly from Equation (\ref{orthh}) that,
for $1\leq i \leq c,$ the
$i$th row of the matrix $K(\omega)$ contains a nonzero entry if and only if $\omega\in S_i,$
i.e., if and only if $i\leq m(\omega).$
And, from Equation (\ref{orthg}), for $c<i\leq c+\widetilde{c}$, the $i$th row
of $K(\omega)$ contains a nonzero entry if and only if $\omega\in \widetilde S_{(i-c)},$
i.e., if and only if $i-c\leq \widetilde m(\omega).$
Therefore, there are exactly $m(\omega)+\widetilde m(\omega)$ nonzero rows in $K(\omega).$
\par
Next, we note that the column indexed by the pair $(j,l)$ has a
nonzero entry only when some $h_{i,j}(\omega_l)$ or $g_{i,j}(\omega_l)$ is nonzero.
That is, the $(j,l)$th column has a nonzero entry only when $\omega_l\in S_j,$ 
i.e., only when $j\leq m(\omega_l).$  So, the maximum number of columns having a nonzero entry in them is
equal to $\sum_{l=0}^{N-1} m(\omega_l),$ which, by the Consistency Equation (\ref{consistency}),
equals $m(\omega)+\widetilde m(\omega),$ and this is exactly the number of rows of $K(\omega)$ that have a nonzero entry.
\par
Therefore, the set of nonzero entries in the matrix $K(\omega)$ are contained in a square submatrix
$L(\omega)$ of size $(m(\omega)+\widetilde m(\omega))\times (m(\omega)+\widetilde m(\omega)).$
\par
Finally, from Equations (\ref{orthg}),(\ref{orthh} and (\ref{orthgh}), we see that the rows of this square matrix $L(\omega)$ are orthonormal.
Hence, the columns of $L(\omega)$ are also orthonormal, and this
implies the orthogonality equations of the theorem.
\end{proof}

\section{Construction of frame wavelets from generalized filters}
We are now ready to use the generalized filters from Section 2 to
extend the construction procedure for frame wavelets described in
\cite{BJ}. Just as in the classical case, the first step of the
construction is to take an infinite product of dilations of the
low-pass filters.

\begin{proposition}
\label{phiprop}
 Let $H=[h_{i,j}]$
 be a generalized low-pass filter as in definition \ref{lowpassdef} Assume that the components of
$H$ are Lipschitz continuous functions in a neighborhood of $0$.
\begin{enumerate}
\item  The infinite product
$$P=\prod_{q=1}^\infty \frac1{\sqrt N}H({B}^{-q}(x))$$
converges almost everywhere on ${\mathbb R}^d,$ and the entries
$P_{i,j}$ of $P$ are square-integrable functions on ${\mathbb
R}^d,$ with $P_{i,j}=0$ for $j>1$.
\item  For $1\leq i \leq c,$ let $\phi_i$ be the inverse
Fourier transform of the function $P_{i,1}.$
Then  the $\widehat{\phi_i}$'s are continuous at 0,
satisfy $\widehat{\phi_{i,1}}(0) = \delta_{i,1},$
and also satisfy the following generalized refinement equation.
$$\widehat{\phi_i}(B(x)) = \frac1{\sqrt N} \sum_{j=1}^c h_{i,j}(x)\widehat{\phi_j}(x).$$
\end{enumerate}
\end{proposition}
\begin{proof} 
Throughout the proof, we will use the following result from linear algebra. If $C$ is a $c\times c$ matrix all of whose eigenvalues have modulus less than $1,$ then $\sum_{k=0}^{\infty}\|C^k(x)\|$ converges, for every $x\in \mathbb R^c.$  This result follows by using the Jordan canonical form for the matrix $C$ to show that $\sum \|C^k(x)\|$ is dominated by a convergent geometric series.
\par We prove the convergence of the infinite product
first. If $\beta$ is the Lipschitz constant, note that by the
low-pass condition we have $\|h_{i,j}(x)\|\leq\beta\|x\|$ for
$\{i,j\}\neq\{1,1\}$, and $\|h_{1,1}(x)-\sqrt N\|\leq\beta \|x\|$,
both for $\|x\|$ sufficiently small. Write $P^k$ for the partial
product $\prod_{q=1}^k \frac1{\sqrt N}H({B}^{-q}(x)),$ where as
before $N= |\det(A)|.$  We first show by induction that for each
fixed $x$, there is a bound $K$ such that $|P_{i,j}^k(x)|\leq K$
for all $1\leq i,j\leq c$, and $k\geq 1.$ To see this, write
\begin{eqnarray*}
|P^k_{i,j}(x)|&=&\sum_{l=1}^cP^{k-1}_{i,l}(x)\frac1{\sqrt N}h_{l,j}({B}^{-k}(x))\\
&\leq& |P_{i,1}^{k-1}(x)|+\left|\sum_{l=2}^cP^{k-1}_{i,l}(x)\right|
\frac{\beta}{\sqrt N} \|{B}^{-k}(x)\|,
\end{eqnarray*}
for $k$ sufficiently large (where in the first term we use the fact that the orthogonality conditions give $|h_{i,1}|\leq
\sqrt N$).  Using induction on $k$, and the linear algebra result mentioned above, we get the bound we seek.
\par
Now, using this bound, we see that for fixed $x$ and for $j\geq
2$,
$|P^k_{i,j}(x)|\leq\left|\sum_{l=1}^cP^{k-1}_{i,l}(x)\right|\frac{\beta}{\sqrt
N} \|{B}^{-k}(x) \|\rightarrow 0$  as $k\rightarrow\infty$.  For
the remaining case of $j=1$, we have
\begin{eqnarray*}
\lefteqn{\left|P^k_{i,1}(x)-P_{i,1}^{k-1}(x)\right|}\\
&=&\left|\sum_{l=1}^cP^{
k-1}_{i,l}(x)\frac1{\sqrt N}h_{l,1}({B}^{-k}x)-P_{i,1}^{k-1}(x)\right|\\
&=&\left|P^{k-1}_{i,1}(x)\left(\frac1{\sqrt
N}h_{1,1}({B}^{-k}x)-1\right)+\sum_{
l=2}^cP^{k-1}_{i,l}(x)\left(\frac1{\sqrt N}h_{l,1}({B}^{-k}x)\right)\right|\\
&\leq&\left|\sum_{l=1}^c|P^{k-1}_{i,l}(x)|\frac{\beta}{\sqrt N}\|{B}^{-k}x\|\right|\\
&<&\frac{\beta}{\sqrt N} cK\|{B}^{-k}x\|.\end{eqnarray*}
We then see that $\{P^k_{i,1}(x)\}$ is Cauchy and thus convergent, 
again by using the linear algebra result mentioned at the beginning of the proof.

To complete the proof of (1), it remains to show that the limiting
functions $P_{i,1}$ are in $L^2(\mathbb R^d)$. To do this, we will
first use induction to prove that
$$\sum_{j=1}^c\int_{{B}^{k}(Q)}|P_{i,j}^k(x)|^2dx\leq 1.$$
We note that because the $h_{k,j}$ are periodic modulo $\mathbb
Z^d,$ the $P_{i,j}^k$ are periodic modulo ${B}^k(\mathbb Z^d),$ and
thus the domain of integration can be replaced by any set that is
congruent to ${B}^{k}(Q)$ modulo ${B}^k (\mathbb Z^d).$ To select
the replacement set, we first choose coset representatives
$\omega_1=0, \omega_2,\cdots, \omega_N$ of $\mathbb
Z^d/{B}(\mathbb Z^d).$ Since
$$\bigsqcup_{n=1}^N(Q+\omega_n)\;\equiv\;{B}(Q)\quad \text{modulo }B(\mathbb Z^d),$$
we can take as our domain of integration the set
${B}^{k-1}(\bigsqcup_{n=1}^N(Q+\omega_n))$ Using this, we have
\begin{eqnarray*}
\lefteqn{\sum_{j=1}^c\int_{{B}^{k}(Q)}|P_{i,j}^k(x)|^2dx}\\
&=&\sum_{j=1}^c\int_{{B}^{k-1}
\bigsqcup_{n=1}^N(Q+\omega_n)}\left(\sum_{l=1}^cP^{k-1}_{i,l}(x)\frac1{\sqrt N}h_{l,j}({B}^{-k}x)\right)\times\\
& &\hskip1.5in\left(\sum_{m=1}^c\overline {P^{k-1}_{i,m}(x)\frac 1{\sqrt N}h_{m,j}({B}^{-l}x)}\right)dx\\
&=&N^{k-1}\int_{\bigsqcup_{n=1}^N(Q+\omega_n)}\sum_{j,l,m}P_{i,l}^{k-1}({B}^{k-1}x)
\overline{P_{i,m}^{k-1}({B}^{k-1}x)}\\
& &\hskip2in\frac1{\sqrt N}h_{l,j}({B}^{-1}x)\overline {\frac1{\sqrt N}h_{m,j}({B}^{-1}x)}dx\\
&=&N^{k-1}\int_{Q}\sum_{l,m}P_{i,l}^{k-1}({B}^{k-1}x)\overline
{P_{i,m}^{k-1}({B}^{k-1}x)}\sum_{j=1}^c\sum_{n=1}^N
\frac1{\sqrt N}h_{l,j}({B}^{-1}x-{B}^{-1}\omega_n)\\
& &\hskip2in \overline {\frac1{\sqrt N}h_{m,j}
({B}^{-1}x-{B}^{-1}\omega_n)}dx.
\end{eqnarray*}

We note now that modulo $\mathbb Z^d,$ the set
$\{{B}^{-1}\omega_1,\;{B}^{-1}\omega_2,\cdots, {B}^{-1}\omega_N\}$
parameterize the distinct $N$ elements of $Q=[-\frac 12, \frac
12)^d$ that map to $(0,0,\cdots, 0)$ under the endomorphism
$\alpha :\mathbb T^d\rightarrow \mathbb T^d,$ where here we are
identifying $\mathbb T^d$ and $[-\frac 12, \frac 12)^d.$ Using the
orthogonality relations satisfied by the $h_{l,j}$, the last
equation simplifies to
\begin{eqnarray*}
N^{k-1}\int_{Q}\sum_{l=1}^cP_{i,l}^{k-1}({B}^{k-1}x)\overline
{P_{i
,l}^{k-1}({B}^{k-1}x)}\chi_{S_l}(x)dx&\leq&\sum_{l=1}^c\int_{{B}^{k-1}(Q)}|P^{k-1}_{i,l}(x)|^2dx\\
&=&\sum_{l=1}^c\int_{{B}^{k-1}(Q)}|P^{k-1}_{i,l}(x)|^2dx.
\end{eqnarray*}
It follows that
$$\sum_{j=1}^c\int_{{B}^{k}(Q)}|P_{i,j}^k(x)|^2dx\leq 1$$ for all $k\in\mathbb N,$ and since
$\cup_{k=0}^{\infty}{B}^{k}(Q)=\mathbb R^d,$ by Fatou's Lemma we
have
$$\sum_{j=1}^c\int_{\mathbb R^d}|P_{i,j}(x)|^2dx\leq 1.$$ In particular we get
$P_{i,j}\in L^2(\mathbb R^d),$ which completes the proof of part
(1).
\par
The refinement equation in part (2) is immediate.
It also follows from the proof above that the infinite product $P$ converges uniformly on neighborhoods of $0,$
and thus that the $\widehat{\phi_i}$'s are continuous at $0.$  Finally, $\widehat{\phi}_i(0)=\delta_{i,1}$ is a consequence
of the low-pass condition.
\end{proof}

We will now use the results of Proposition \ref{phiprop} to build
a frame wavelet. We begin by generalizing the computational ideas
in \cite{BJ}. Define two Hilbert spaces ${\mathcal{H}}$ and
$\widetilde{{\mathcal H}}$ by ${\mathcal{H}} = \bigoplus_{j=1}^c L^2(S_i)$
and $\widetilde{{\mathcal{H}}} = \bigoplus_{k=1}^{\widetilde c}
L^2(\widetilde S_k);$ and two operators $S_H:{\mathcal{H}}\to {\mathcal{H}}$
and $S_G:\widetilde{{\mathcal{H}}}\to {\mathcal{H}}$ by:
$$[S_H(f)](\omega) = H^t(\omega)f(\alpha(\omega))$$
and
$$[S_G(\widetilde f)](\omega) = G^t(\omega)\widetilde f(\alpha(\omega)),$$
where $f\equiv \bigoplus f_j \in {\mathcal{H}},$ $H^t$ and $G^t$ denote
the transposes of the matrix functions $H$ and $G$ respectively,
and as above $\alpha$ denotes the map on the torus ${\mathbb T}^d$
induced by the action of the transpose $B$ of $A$ acting on
${\mathbb R}^d.$ Note that $f_j(\alpha(\omega))$ and $\widetilde
f_k(\alpha(\omega))$ are necessarily only defined when
$\alpha(\omega)$ belongs to $S_j$ for the first case and
$\widetilde S_k$ for the second.  We define these functions to be
$0$ outside of these domains. We record here the formulas for the
adjoints of the two operators:
\begin{equation}
[S_H^*(f)](\omega) = \frac1N \sum_{l=0}^{N-1} \overline H(\omega_l)f(\omega_l),\label{SH*}
\end{equation}
and
$$[S_G^*(f)](\omega) = \frac1N \sum_{l=0}^{N-1} \overline G(\omega_l)f(\omega_l).$$
It will also be convenient to have explicit formulas for the powers of both $S_H$ and $S_H^*:$
$$[S_H^n(f)](\omega) = \prod_{k=0}^{n-1} H^t(\alpha^k(\omega)) f(\alpha^n(\omega),$$
and
\begin{equation}
[{S_H^*}^n(f)](\omega) = \frac1{N^n} \sum_{l=0}^{N^n-1} \prod_{ k=n-1}^0 \overline H(\alpha^k(\omega_{l,n}))f(\omega_{l,n}),\label{SH*n}
\end{equation}
where the $\omega_{l,n}$'s are the $N^n$ points $\zeta\in {\mathbb
T}^d$ for which $\omega=\alpha^n(\zeta).$
\par
We also include here an estimate involving these operators that we
will need later, \begin{equation}\sum_{i=1}^c \sum_{j=1}^c
\sum_{l= 0}^{N^n-1}\left|\left[\prod_{k=0}^{n-1}
H^t(\alpha^k(\omega_{l,n}))\right]_{i,j}\right|^2 \leq
cN^n,\label{estimate} \end{equation}
 which we prove using
induction. The case $n=1$,
$$ \sum_{i=1}^c \sum_{j=1}^c \sum_{l= 0}^{N-1} |h_{j,i}(\omega_{l})|^2 \leq cN,$$
follows from equation (\ref{orthh}) together with the fact that
$\sum_{j=1}^c\chi_{S_j}(\omega)=m(\omega)\leq c.$  For the
induction step, note that the $\omega_{l,n+1}$ can be labeled in
such a way that $\alpha^n(\omega_{sN+q,n+1})=\omega_q$, so that
$\omega_{sN+q,n+1}=(\omega_q)_{s,N}.$ Thus, writing $l=sN+q$ and
using the Cauchy-Schwarz inequality, we obtain:
\begin{eqnarray*} \lefteqn{\sum_{i=1}^c \sum_{j=1}^c \sum_{l= 0}^{N^{n+1}-1}\left|\left[\prod_{k=0}^{n} H^t(\alpha^k(\omega_{l,n+1}))\right]_{i,j}\right|^2}\\
&=&\sum_{i=1}^c \sum_{j=1}^c
\sum_{l=0}^{N^{n+1}-1}\left|\sum_{r=1}^c\left[\prod_{k=0}^{n-1}
H^t(\alpha^k(\omega_{l,n+1}))\right]_{i,r}\left[H^t(\omega_q)\right]_{r,j}\right|^2\\
&\leq&\sum_{i=1}^c \sum_{j=1}^c \sum_{s=0}^{N^{n}-1}\sum_{q=0}^{N-1}\left(\sum_{r=1}^c\left|\prod_{k=0}^{n-1}h_{r,i}(\alpha^k(\omega_{sN+q,n+1})\right|^2\right)\left(\sum_{r'=1}^c\left|h_{j,r'}(\omega_q)\right|^2\right)\\
&\leq& cN^{n} N
\end{eqnarray*}

\begin{lemma}
\label{Cuntz}
The operators $S_H$ and $S_G$ satisfy the
following relations: \begin{enumerate}
\item\hskip0.5em $S_H^*S_H = I,$ the identity operator on ${\mathcal{H}}.$
\item\hskip0.5em $S_G^*S_G = \widetilde I,$ the identity operator on $\widetilde{{\mathcal{H}}}.$
\item\hskip0.5em $S_H^*S_G = 0$ and $S_G^*S_H = \widetilde{0},$ the $0$ operators on the appropriate Hilbert spaces.
\item\hskip0.5em $S_HS_H^* + S_GS_G^* = I,$ the identity operator on ${\mathcal{H}}.$
\end{enumerate}
\end{lemma}
\begin{proof} These are direct consequences of the orthogonality
relations satisfied by generalized filter systems relative to $m$
and $\widetilde m.$ For example, to prove (1), for $f\in\mathcal{H}$,
we write
\begin{eqnarray*}
S_H^*S_Hf(\omega)&=& \frac1N\sum_{l=0}^{N-1}
\overline{H(\omega_l}S_H f(\omega_l)\\
&=&\frac1N\sum_{l=0}^{N-1}\overline{H(\omega_l)}H^t(\omega_l)f(\omega)\\
&=&f(\omega), \end{eqnarray*} where the last equality follows
from \ref{orthh}.  The other parts of the lemma follow similarly
from \ref{orthg}, \ref{orthgh}, and \ref{orthcol} respectively.
\end{proof}

We note   that both $S_H$ and $S_G$ are partial isometries,
and that the relations in the lemma are similar to the famous
relations defining the Cuntz algebra ${\mathcal{O}}_2$.
We now use these operators to decompose the Hilbert space ${\mathcal{H}}$
in a convenient way.

 \begin{lemma}
\label{calHlemma}
 Let $R_0$ denote the
range of $S_G,$ and let $R_n=S_H^n(R_0).$ Then ${\mathcal{H}} =
\bigoplus_{n=0}^\infty  R_n.$ Moreover, if $z$ is any element in
${\mathbb Z}^d,$ and $\widetilde f_{k,z}$ is the element of
$\widetilde{{\mathcal{H}}}$ whose $k$th component is the exponential
function $e^{2\pi i\langle \omega\mid z\rangle}$ and whose other
components are 0, then the collection $\{S_H^n(S_G(\widetilde
f_{k,z}))\},$ for $k$ running from 1 to $\widetilde c$ and each
$z$ running through ${\mathbb Z}^d,$ forms a Parseval frame for
the subspace $R_n.$ Consequently, the collection
$\{S_H^n(S_G(\widetilde f_{k,z}))\},$ $1\leq k \leq \widetilde c,$
$z\in {\mathbb Z}^d,$ and $0\leq n <\infty,$ forms a Parseval
frame for ${\mathcal{H}}.$
\end{lemma}
\begin{proof} That the subspaces $\{R_n\}$ are orthogonal follows
directly from the relations in Lemma 3.2. That the elements
$\{S_H^n(S_g(\widetilde f_{k,z}))\}$ form a Parseval frame for
$R_n$ follows from the fact that $S_H$ and $S_G$ are partial
isometries, together with the fact that the elements $\{\widetilde
f_{k,z}\},$ as $k$ runs from 1 to $\widetilde c$ and $z$ runs over
${\mathbb Z}^d,$ form a Parseval frame for $\widetilde{{\mathcal{H}}}.$
\par
Write $R_\infty$ for the orthogonal complement in ${\mathcal{H}}$ of the direct sum $\bigoplus R_n.$
We must show that $R_\infty = \{0\}.$
Note that $R_\infty$ is invariant under both $S_H$ and $S_H^*,$ and therefore that the restriction
of both these operators to $R_\infty$ are unitary operators.
\par
By way of contradiction, suppose that $f_0$ is a unit vector in
$R_\infty.$ For each natural number $n,$ write
$f_n={S_H^*}^n(f_0).$ Note that $f_{n+m} = {S_H^*}^m(f_n).$ We
make two observations about $f_n.$ First of all, for almost all
$\omega\in {\mathbb T}^d,$  we have
\begin{eqnarray*}
\|f_n(\omega)\|^2&=&\|{S_H^*}^n(f_0(\omega))\|^2 \cr
&=&\frac1{N^{2n}} \sum_{i=1}^c \left|\sum_{p=1}^c
\sum_{l=0}^{N^n-1}\left[\prod_{k=n-1}^0 \overline
H(\alpha^k(\omega_{l,n}))\right]_{i,p}
[{f_0}]_p(\omega_{l,n})\right|^2 \cr &\leq& \frac1{N^n}
\sum_{i=1}^c \sum_{p=1}^c \sum_{l=0}^{N^n-1}
\left|\left[\prod_{k=n-1}^0 \overline
H(\alpha^k(\omega_{l,n}))\right]_{i,p}\right|^2 \times \frac1{N^n}
\sum_{p'=1}^c \sum_{l'=0}^{N^n-1}
\left|[{f_0}]_{p'}(\omega_{l',n})\right|^2 \cr &\leq& c\times
\frac1{N^n} \sum_{p'=1}^c \sum_{l'=0}^{N^n-1}
|[{f_0}]_{p'}(\omega_{l',n})|^2,\end{eqnarray*} where the last
inequality follows from the transpose of equation \ref{estimate}.  By the
pointwise ergodic theorem this then implies that
$$\limsup \|f_n(\omega)\|^2 \leq c\int_0^1 \|f_0(\omega)\|^2\, d\omega = c.$$
\par
Secondly, since $S_H$ is unitary on $R_\infty,$ we have
\begin{eqnarray*} \frac1{N^n} \sum_{l=0}^{N^n-1} \|f_0(\omega_{l,n})\|^2
&=& \frac1{N^n} \sum_{l=0}^{N^n-1}
\|S_H^n{S_H^*}^n(f_0)(\omega_{l,n})\|^2 \cr &=& \frac1{N^n}
\sum_{l=0}^{N^n-1} \sum_{i=1}^c\left|\sum_{p=1}^c
\left[\prod_{k=0}^{n-1} H^t(\alpha^k(\omega_{l,n}))\right]_{i,p}
[{S_H^*}^n(f_0)]_p(\omega)\right|^2 \cr &\leq& \frac1{N^n}
\sum_{i=1}^c \sum_{l=0}^{N^n-1}\sum_{p=1}^c
\left|\left[\prod_{k=0}^{n-1}
H^t(\alpha^k(\omega_{l,n}))\right]_{i,p}\right|^2 \times
\sum_{p'=1}^c |[{f_n}]_{p'}(\omega)|^2 \cr &\leq& c
\|f_n(\omega)\|^2,\end{eqnarray*} implying (again by the pointwise
ergodic theorem) that
$$\liminf \|f_n(\omega)\|^2 \geq \frac 1c.$$
Consequently, by Egorov's Theorem, for any $\epsilon>0,$ there
exists an $M_0$ and a set $E\subseteq {\mathbb T}^d$ such that the
measure of the complement of $E$ is less than $\epsilon,$ and
$\frac 1c-\epsilon < \|f_n(\omega)\|^2 < c+\epsilon$ for all
$\omega\in E$ and all $n\geq M_0,$.  Before we apply this theorem,
we will establish some other estimates needed in our choice of
$\epsilon$.

First, we pick an integer $K_0\geq 3 \log_2c+9$, so that we then
have
\begin{equation} \label{K0condition} \sqrt N^{K_0+1}>32c^{\frac32}.
\end{equation}(This follows since $N=|\det
A|\geq 2$.)

Next, we choose a $\delta>0$ so that
\begin{equation}
\label{delta1}
\delta< \frac 1{4c^{K_0+3}\sqrt N^{K_0}},
\end{equation}
and
\begin{equation}
\label{delta2}
\sqrt N(1-\delta)^{K_0+1} -\delta c^{K_0}\; > \frac{\sqrt N}2.
\end{equation}
Note the second condition is possible since the function on the
left hand side approaches $\sqrt{N}$ as $\delta$ goes to $0.$
\par
From the low-pass condition and the requirement of Lipschitz near
0, we know that for any $\delta>0$ there exists a neighborhood $U$
of $0\in {\mathbb T}^d$ such that $|h_{i,j}(\omega)|<\delta$ for
all $\omega\in U$ and all pairs $(i,j)\neq (1,1),$ and
$|h_{1,1}(\omega)|> \sqrt N(1-\delta)$ for all $\omega\in U.$ Let
$U$ be the neighborhood corresponding to our choice of $\delta$
satisfying (1) and (2) above.  By continuity of $\alpha,$ there
further exists a neighborhood $V\subseteq U$ such that for every
$\omega\in V$ we have $\alpha^k(\omega)\in U$ for all $0\leq k
\leq 2(K_0+1).$ Hence, if $\widetilde P(\omega)$ is the matrix
given by $\widetilde P(\omega)=\prod_{k=0}^{5}
H^t(\alpha^k(\omega)),$ then for $(i,j)\neq (1,1),$ $0\leq l\leq
K_0+1$ and $\omega\in V$,
$$|\widetilde P_{i,j}(\alpha^l(\omega))| \leq \delta (c\sqrt
N)^{K_0}.$$ (There are $c^{K_0}$ summands, each having $K_0$
factors, and in each, one factor is bounded by $\delta$, and the
other factors (by \ref{orthh})
 are bounded by $\sqrt N.$) For $(i,j)=(1,1),$ we
have
$$(\sqrt N(1-\delta))^{K_0+1} -\delta(c\sqrt N)^{K_0} \leq |\widetilde P_{1,1}(\alpha^l(\omega))| \leq \sqrt N^{K_0+1} + \delta(c\sqrt N)^{K_0}.$$
(Again, there are $c^{K_0}$ summands, each having $K_0+1$ factors.
One of these summands is bounded below by $ (\sqrt
N(1-\delta))^{K_0+1}$ and above by $\sqrt N^{K_0+1},$ and the
other $(c^{K_0}-1$ summands are bounded by $\delta\sqrt N^{K_0}.$)
\par
Now choose an  $\epsilon$ smaller than the measure of $V,$ and
also smaller than $1/(4c),$ so that we are assured that the
corresponding set $E$ will satisfy $\lambda(E\cap V) >0.$ In fact,
we may even assume that the set of $\omega \in V$ for which
$\alpha^k(\omega)\in E\cap U$ for all $0\leq k \leq 2(K_0+1)$ has
positive measure.
\par
Fix an $\omega\in E\cap V$ for which $\alpha^k(\omega)\in E\cap U$
for $0\leq k \leq 2(K_0+1).$ If $M_0$ is the natural number
corresponding to this choice of $\epsilon$ as above, and $n\geq
M_0,$ we have the following estimates on the components of
$f_n(\omega).$ First, for $i\neq 1$ and $0\leq l\leq K_0+1,$
\begin{eqnarray*} |{f_n}_i(\alpha^l(\omega))| &=& |[S_H^{K_0+1}(f_{n+K_0+1})]_i(\alpha^l(\omega))| \cr
&=& \left|\sum_{j=1}^c \widetilde
P_{i,j}(\alpha^l(\omega))[{f_{n+K_0+1}}]_j(\alpha^{l+K_0+1}(\omega))\right|
\cr &\leq& \sum_{j=1}^c \left|\widetilde P_{i,j}(\alpha^l(\omega))\right|
\times\left|[{f_{n+K_0+1}}]_j(\alpha^{l+K_0+1}(\omega))\right| \cr &\leq&
c\delta(c\sqrt N)^{K_0}\sqrt{c+\epsilon} \cr &\leq& 2\delta
c^{K_0+2} \sqrt N^{K_0},\end{eqnarray*} so that by condition (\ref{delta1})
on $\delta$,

$$\sum_{i=2}^c |[{f_n}]_i(\alpha^l(\omega))|^2 \leq 4c\delta^2 c^{2(K_0+2)} N^{K_0} < \frac 1{2c}.$$
Therefore, because $\omega\in E,$ for all $n\geq N_0$ we must
have, for $0\leq l\leq K_0+1$,
$$|[{f_n}]_1(\alpha^l(\omega))|^2 > \frac1{2c} - \epsilon > \frac1{4c}.$$
\par
On the other hand,
\begin{eqnarray*} [{f_n}]_1(\omega) &=& \sum_{j=1}^c \widetilde P_{1,j}(\omega) [{f_{n+K_0+1}}]_j(\alpha^{K_0+1}(\omega)) \cr
&=& \widetilde
P_{1,1}(\omega)[{f_{n+K_0+1}}]_1(\alpha^{K_0+1}(\omega)) +
\sum_{j=2}^c \widetilde
P_{1,j}(\omega)[{f_{n+K_0+1}}]_j(\alpha^{K_0+1}(\omega))
,\end{eqnarray*} implying that
\begin{eqnarray*}|[{f_{n+K_0+1}}]_1(\alpha^{K_0+1}(\omega))|
&\leq& \frac{|[{f_n}]_1(\omega)| + \sum_{j=2}^c \left|\widetilde
P_{1,j}(\omega)\right|\times
\left|[f_{n+K_0+1}]_j(\alpha^{K_0+1}(\omega))\right|}{|\widetilde
P_{1,1}(\omega)|} \cr &\leq& \frac{\sqrt{c+\epsilon}+c\delta(c\sqrt N)^{K_0}\sqrt{c+\epsilon}}{(\sqrt
N(1-\delta))^{K_0+1} -\delta(c\sqrt N)^{K_0}}\cr &\leq&
\frac{4\sqrt{c+\epsilon}}{\sqrt N^{K_0+1}}\cr
&\leq&\frac{8\sqrt c}{\sqrt N^{K_0+1}}\cr &\leq& \frac1{4c},
\end{eqnarray*}
where in the third step we use conditions (\ref{delta1}) and (\ref{delta2}) to simplify
the numerator and denominator respectively, and in the final step, we use condition (\ref{K0condition}). Hence, the
point $\alpha^{K_0+1}(\omega)$ satisfies
$$\frac1{4c} \leq |[{f_{n+K_0+1}}]_1(\alpha^{K_0+1}(\omega))| ^2 < \frac 1{4c},$$
which is a contradiction.
\end{proof}

We now state our main result.
\begin{theorem}
\label{waveletGMRAth}
Let $\phi_i$ be defined from the infinite
product of the low-pass filter system $H$ as in Proposition
\ref{phiprop}, and let $G=[g_{k,j}]$ be the corresponding high-%
pass filter system. Then, if we define a function $\psi_k\in L^2({\mathbb R}^d)$,
for $1\leq k \leq \widetilde c,$  by
$$\widehat{\psi_k}(x) = \frac1{\sqrt N}\sum_{j=1}^c g_{k,j}({B}^{-1}(x))\widehat{\phi_j}({B}^{-1}(x)),$$
the collection $\{\psi_k\}$ is a Parseval frame wavelet for
$L^2({\mathbb R}^d)$ relative to dilation by $A$ and translation
by lattice points $z\in{\mathbb Z}^d.$  Further, if $V_0$ is the closed linear span of
the translates of the $\phi_i$'s by lattice elements $z\in{\mathbb
Z}^d,$ then $\{V_j\}\equiv\{\delta^j(V_0)\}_{j\in\mathbb Z}$ is a
generalized multiresolution analysis for $L^2({\mathbb R}^d)$.
\end{theorem}

\begin{proof}
We prove first that the $\{\psi_k\}$ form a frame wavelet. For
convenience in what follows, we introduce the following notation.
For each $1\leq k\leq \widetilde c$ and each $z\in{\mathbb Z}^d,$
write $\gamma^{k,z}$ for the element $S_G(\widetilde f_{k,z})$ of
$R_0\subseteq {\mathcal{H}},$ where as in the previous lemma,
$\widetilde f_{k,z}$ is the element of $\widetilde{{\mathcal{H}}}$
whose $k$th component is the exponential function $e^{2\pi
i\langle \omega\mid z\rangle}$ and whose other components are 0.
  Note that $\gamma^{k,z} = \bigoplus_{j=1}^c
\gamma^{k,z}_j,$ where
$$\gamma^{k,z}_j(\omega) = g_{k,j}(\omega)e^{2\pi i\langle \alpha(\omega)\mid z\rangle}.$$
It then follows from Lemma 3.3 that for each $F\in{\mathcal{H}}$ we have
$$\|F\|^2 = \sum_{n=0}^\infty   \sum_{k=1}^{\widetilde c} \sum_{z\in {\mathbb Z}^d}
 |\langle F\mid S_H^n(\gamma^{k,z})\rangle|^2.$$
\par
For each integer $n,$ each $1\leq k \leq\widetilde c,$ and each
$z\in{\mathbb Z}^d,$ define the function $\psi_{n,k,z}$ by
$$\psi_{n,k,z}(t) = \sqrt N^n \psi_k(A^n(t)+z),$$
and note that the Fourier transform of $\psi_{n,k,z}$ is given by
$$\widehat{\psi_{n,k,z}}(x)
=N^{\frac{-n-1}2} e^{2\pi i\langle {B}^{-n}(x )\mid z\rangle}
\sum_{j=1}^c g_{k,j}({B}^{-n-1}(x )) \widehat{\phi_j}({B}^{-n-1}(x
)).$$ We wish to prove that the collection $\{\psi_{n,k,z}\}$ is a
Parseval frame for $L^2({\mathbb R}^d).$
\par
Now, let $f$ be an element of $L^2({\mathbb R}^d)$ whose Fourier
transform $\widehat f$ has compact support. For a fixed integer
$J\geq 0,$ we have
\begin{eqnarray*} \lefteqn{\sum_{n=-\infty}^J\sum_{k=1}^{\widetilde c} \sum_{z\in {\mathbb Z}^d}
|\langle f\mid \psi_{n,k,z}\rangle|^2}\cr &=&  \sum_{n=-\infty}^J
\sum_{k=1}^{\widetilde c} \sum_{z\in{\mathbb Z}^d} |\langle
\widehat f \mid \widehat{\psi_{n,k,z}}\rangle|^2 \cr &=&
\sum_{n=-\infty}^J \sum_{k=1}^{\widetilde c} \sum_{z\in {\mathbb
Z}^d} \left|\int_{{\mathbb R}^d} \widehat{f}(x)
\overline{N^{\frac{-n-1}2} \sum_{j=1}^c e^{2\pi i\langle
{B}^{-n}(x)\mid z\rangle} g_{k,j}({B}^{-n-1}(x))
\widehat{\phi_j}({B}^{-n-1}(x))}\, dx \right|^2 \cr &=&
\sum_{n=-\infty}^{J} \sum_{k=1}^{\widetilde c} \sum_{z\in {\mathbb
Z}^d} N^{-n+2J+1} \left|\int_{{\mathbb R}^d}
\widehat{f}({B}^{J+1}(x)) \overline{e^{2\pi i\langle
{B}^{-n+J+1}(x) \mid z \rangle}}\right.\times\cr
&&\hskip2in\left.\overline{\sum_{j=1}^c
g_{k,j}({B}^{-n+J}(x)) \widehat{\phi_j}({B}^{-n+J}(x)) }\, dx
\right|^2 \cr &=& \sum_{n=0}^\infty \sum_{k=1}^{\widetilde c}
\sum_{z\in {\mathbb Z}^d} N^{n+1+J}\left|\int_{{\mathbb R}^d}
\widehat{f}({B}^{J+1}(x)) \overline{e^{2\pi i\langle {B}^{n+1}(x)
\mid z \rangle} \sum_{j=1}^c g_{k,j}({B}^n(x))
\widehat{\phi_j}({B}^n(x)) }\, dx \right|^2 \cr &=&
\sum_{n=0}^\infty \sum_{k=1}^{\widetilde c} \sum_{z\in {\mathbb
Z}^d} N^{1+J}\left |\int_{{\mathbb R}^d} \widehat{f}({B}^{J+1}(x)
\overline{e^{2\pi i\langle {B}^{n+1}(x) \mid z\rangle}
\sum_{j=1}^c g_{k,j}({B}^n(x))}\right.\times\cr
&&\hskip3in\left.\overline{ \left[ \prod_{p=n-1}^0
H({B}^p(x)\Phi(x)\right]_j}\, dx \right|^2 \cr &=&
\sum_{n=0}^\infty \sum_{k=1}^{\widetilde c} \sum_{z\in {\mathbb
Z}^d} N^{1+J} |\int_{{\mathbb T}^d} \sum_{j=1}^c \sum_{\zeta\in
{\mathbb Z}^d} \widehat{f}({B}^{J+1}(x+\zeta))
\overline{\widehat{\phi_j}(x+\zeta)} \overline{ \prod_{p=0}^{n-1}
H^t({B}^p(x)} \overline{\gamma^{k,z}}({B}^n(x))\, dx|^2 \cr &=&
\sum_{n=0}^\infty \sum_{k=1}^{\widetilde c} \sum_{z\in {\mathbb
Z}^d} |\langle F^J\mid S_H^n(\gamma^{k,z})\rangle_{{\mathcal{H}}}|^2
\cr &=& \|F^J\|^2_{{\mathcal{H}}},\end{eqnarray*} where
$F^J=\bigoplus {F^J}_j$ is the element of ${\mathcal{H}}$ given by
$${F^J}_j(\omega) = \sqrt N^{1+J} \sum_{\zeta\in{\mathbb Z}^d} \widehat f({B}^{1+J}(x+\zeta)) \overline{\widehat{\phi_j}}(x+\zeta).$$
Hence,
\begin{eqnarray*}\lefteqn{ \sum_{n=-\infty}^J \sum_{k=1}^{\widetilde c} \sum_{z\in{\mathbb Z}^d}
|\langle f\mid \psi_{n,k,z}\rangle|^2}\cr &=& \sum_{j=1}^c
\int_{{\mathbb T}^d} |{F^J}_j(\omega)|^2\, d\omega \cr &=&
\sum_{j=1}^c \int_{{B}^{1+J}({\mathbb T}^d)} |\sum_{\zeta\in
{\mathbb Z}^d} \widehat f(x+{B}^{1+J}(\zeta))
\overline{\widehat{\phi_j}}({B}^{-1-J}(x)+\zeta)|^2\, dx  \cr &=&
\sum_{j=1}^c \int_{{\mathbb R}^d} \chi_{{B}^{1+J}({\mathbb
T}^d)}(x) \left|\sum_{\zeta\in {\mathbb Z}^d} \widehat
f(x+{B}^{1+J}(\zeta)) \overline{\phi_j}({B}^{-1-J}(x)+\zeta)\right|^2\,
dx .\end{eqnarray*} Now, because $\widehat f$ has compact
support, and the matrix $B$ is expansive, there exists a $J'$ such
that the support of $\widehat f$ is contained in ${B}^{1+J'}(Q).$
There must also exist, by the compactness of
${B}^{1+J'}\overline{(Q)}$, a $J_0$ such that ${B}^{1+J}(Q)$
contains ${B}^{1+J'}(Q)$ for all $J\geq J_0.$ Now, for $J\geq
J_0$, the product $\chi_{{B}^{1+J}(Q)}(x) \widehat
f(x+{B}^{1+J}(\zeta))$ is nonzero only if $x={B}^{1+J}(y)$ for
some $y\in Q$, and also $x+{B}^{1+J}(\zeta)={B}^{1+J}(y+\zeta)\in
{B}^{1+J}(Q).$ Consequently, for any $J>J_0,$ and any $\zeta\neq
0,$ we must have $\chi_{{B}^{1+J}(Q)}(x) \widehat
f(x+{B}^{1+J}(\zeta)) = 0$ for all $x\neq 0$ Hence, for $J>J_0$ we
have
$$ \sum_{n=-\infty}^J \sum_{k=1}^{\widetilde c} \sum_{z\in {\mathbb Z}^d}
|\langle f\mid \psi_{n,k,z}\rangle|^2 = \sum_{j=1}^c
\int_{{\mathbb R}^d} |\widehat f(x)|^2
|\widehat{\phi_j}({B}^{-1-J}(x))|^2\, dx,$$ so
that, by the Dominated Convergence Theorem, we obtain
\begin{eqnarray*} \sum_{n=-\infty}^\infty \sum_{k=1}^{\widetilde c} \sum_{z\in{\mathbb Z}^d}
|\langle f\mid \psi_{n,k,z}\rangle|^2 &=& \int_{{\mathbb R}^d}
|\widehat f(x)|^2 \sum_{j=1}^c |\widehat{\phi_j}(0)|^2\, dx \cr
&=& \|f\|^2.\end{eqnarray*} This demonstrates the Parseval frame
property for elements $f$ whose Fourier transforms have compact
support. For general $f\in L^2({\mathbb R}^d),$ the Parseval
equality follows from the density of these functions.

It remains to show that $\{V_j\}$ form a GMRA.  Properties (2) and
$(4')$ of Definition \ref{GMRAdef} are direct consequences of the definition of the $V_j,$ and
property (1) follows immediately from Proposition \ref{phiprop} (2).   The
density required in Property (3) is a consequence of the fact that
the $\{\psi_{n,k,z}\}$ form a Parseval frame for $L^2(\mathbb
R^d).$

To finish the proof, we now show that $\cap V_j=\{0\}$.
Write $\{\phi'_i\}$ for a sequence of elements whose translates
form a Parseval frame for $V_0$ (such a sequence must exist since
$V_0$ is closed under translation), and $P_j$ for the orthogonal
projection operator onto the subspace $V_j.$
To prove that  $\cap V_j=\{0\}$, it will suffice to show
that $\lim_{j\to\infty}\|P_{-j}(f)\|=0$
for each $f\in L^2(\mathbb R^d).$
By a standard approximation argument, it will suffice to
show this holds on a
dense subset of $L^2(\mathbb R^d).$
Thus, let $f$ be a Schwartz function for which $\widehat f$
vanishes in some neighborhood $N_f$ of 0,
and write $C_f$ for the (finite) number $\sum_k|f*f^{*}(k)|.$
Such $f$'s are dense in $L^2(\mathbb R^d).$
The Poisson Summation Formula holds for such an $f,$ and we will use it
in the following form:
$$\frac 1{N^j}\sum_l|\widehat f(B^{-j}(\xi + l))|^2
=\sum_ kf*f^{*}(A^j(k))e^{-2\pi i\langle k\mid\xi\rangle}.$$
Now, for each $\xi\in {\mathbb T}^d,$ let $l_j(\xi )$
be the largest number for which $B^{-j}(\xi + l)\in N_f$ for
all $\|l\|<l_j(\xi ).$
Because $A$ (and thus $B=A^t$) is expansive, we must have that $l_j(\xi )$ tends to
infinity for almost every $\xi .$
Finally, we use the fact  that the function
$$m(\omega )=\sum_i\chi_{S_i}(\omega )=\sum_i\sum_l|\widehat{\phi_
i'}(\omega + l)|^2$$
is  integrable on the cube.
Hence, we have
\begin{eqnarray*}
\|P_{-j}(f)\|^2
&=& N^j\sum_i\sum_{ z}\left| \int_{\mathbb R^d}\widehat{\phi_i'}(B^j(\xi ))
\overline {\widehat f}(\xi )e^{-2\pi i\langle B^j(\xi )\mid z\rangle}\,d\xi \right| ^2\cr
&=&\frac1{N^j}  \sum_i\sum_{z}\left| \int_{\mathbb R^d}\widehat{
\phi_i'}(\xi )
 \overline {\widehat f}(B^{-j}(\xi ))e^{-2\pi i\langle z\mid\xi\rangle}\,d\xi \right| ^2\cr
&=&\frac 1{N^j}\sum_i\int_{ {\mathbb T}^d}\left| \sum_l\widehat{\phi_i'}(\xi + l)
 \overline {\widehat f}(B^{-j}(\xi + l))\right| ^2\,d\xi\cr
&=&\frac 1{N^j}\sum_ i\int_{{\mathbb T}^d}\left| \sum_{\| l\| \geq l_j(\xi )}\widehat{\phi_i'}(\xi
+ l)
 \overline {\widehat f}(B^{-j}(\xi + l))\right| ^2\,d\xi\cr
&\leq& \frac 1{N^j}\sum_i\int_{{\mathbb T}^d}\sum_{|l|\geq l_
j(\xi )}|\widehat{\phi_i'}(\xi + l)|^2
\sum_{|l|\geq l_j(\xi
)}|\widehat f(B^{-j}(\xi + l))|^2\,d\xi\cr
&=& \int_{{\mathbb T}^d}\left[\sum_i\sum_{|l|\geq l_j(\xi )}|\widehat{\phi_
i'}(\xi + l)|^2\right]
\left[\frac 1{N^j}\sum_l|\widehat f(B^{
-j}(\xi + l))|^2\right]\,d\xi\cr
&=& \int_{{\mathbb T}
^d}\left[\sum_i\sum_{|l|\geq l_j(\xi )}|\widehat{\phi_i'}(\xi + l)|^2\right]
\left[\sum_kf*f^{*}(A^j(k))e^{-2\pi i\langle k\mid\xi\rangle}\right]\,d\xi\cr
&\leq& \int_{{\mathbb T}^d}\left[\sum_i\sum_{|l|\geq l_
j(\xi )}|\widehat{\phi_i'}(\xi + l)|^2\right]\left[\sum_k|f*f^{*}(
A^j(k))|\right]\,d\xi\cr
&\leq& C_f\int_{{\mathbb T}^d}\sum_i\sum_{|l|\geq
l_j(\xi )}|\widehat{\phi_i'}(\xi + l)|^2\,d\xi ,\end{eqnarray*}
which approaches 0 as $j$ goes to infinity by the
Dominated Convergence Theorem,
the integrand here being bounded by $m.$
\end{proof}

Theorem \ref{waveletGMRAth} shows that the functions $\psi_k$ are
in the set $V_1$, so that the Parseval frame wavelet
$\{\psi_{k}\}$ that we have constructed is \textit{obtained from}
the GMRA $\{V_j\}$, in the sense defined by Zalik \cite{Zl}. It
is an open question whether the wavelet is \textit{associated
with} the GMRA $\{V_j\}$ in the sense that $V_j$ is the closure of
the span of $\{\delta^l(\tau_z(\psi_{k}))\}_{l<j}$ (see e.g.
\cite{Bw}). Moreover, the multiplicity function of the GMRA
produced by the theorem may not coincide with the multiplicity
function used in the construction of the filters, as the first
of the examples below (Example \ref{ExMult1}) shows.

\begin{example}
\label{ExMult1}Let $d=1$, $A=2$, and the multiplicity functions $m$ and
$\widetilde{m}$ both be identically 1, so that we start in the
classical setting.  For our filters, we take
$$h=\sqrt2\left[\chi_{-\frac18,\frac18)} + \chi_{\pm[\frac14,\frac38)}\right],$$  and
$$g=\sqrt2\left[\chi_{\pm\frac18,\frac14)} + \chi_{\pm[\frac38,\frac12)}\right].$$
Note that $h$ and $g$ satisfy the definitions of generalized low-
and high-pass filters, and are Lipschitz in a neighborhood of 0.
The infinite product $P=\prod_{j=1}^\infty \frac1{\sqrt2} h(\frac
x{2^j})$ equals $\chi_{[-\frac14,\frac14)}.$ The integer
translates of the function $\phi$ whose Fourier transform is this
infinite product $P,$ are not orthonormal, and do not determine
the core subspace of any classical multiresolution analysis. On
the other hand, as predicted by Theorem \ref{waveletGMRAth}, the
standard construction in this case does produce a generalized
multiresolution analysis $\{V_j\},$ with multiplicity function
$m=\chi_{[-\frac 14,\frac14)}$.  As guaranteed by Theorem
\ref{waveletGMRAth}, the construction also produces a Parseval frame
wavelet $\psi,$ here given by
$$\widehat{\psi}(x) = \chi_{[-\frac12,-\frac14)} + \chi_{[\frac14,\frac12)}.$$
In this case the function $\phi$ is easily constructed out of
negative dilates of $\psi$, so the wavelet is necessarily
associated with as well as obtained by the GMRA.
\end{example}

We end this section with an another example for dilation by 2 in $L^2(\mathbb R)$.  This one begins 
with filters for the multiplicity function of the well-known Journ\'e wavelet, whose
Fourier transform is 
the characteristic function of the set 
$$[-\frac{16}{7},-2)\cup [-\frac{1}{2},-\frac{2}{7})\cup [\frac{2}{7},\frac
{1}{2}]\cup [2,\frac{16}{7}).$$  We use Theorem \ref{waveletGMRAth} to build a Parseval
wavelet with a $C^{\infty}$ Fourier transform that is associated with the Journ\'e multiplicity function.  

\begin{example}
\label{ExJourne}The Journ\'e multiplicity function $m$ takes on the values $0,1,$ and $2,$ 
with
$S_1=[-\frac {1}{2},-\frac{3}{7})\cup [-\frac {2}{7},\frac {2}{7})\cup 
[\frac {3}{7},\frac{1}{2})$ and $S_2=[-\frac {1}{7},\frac{1}{7}).$ 
The complementary 
multiplicity function $\widetilde{m}(x)\equiv 1,$ since the Journ\'e wavelet is a single orthonormal wavelet.

To build the filters, we let $p_0$ be a $C^\infty$ classical (MRA) low-pass filter for dilation by 2 (\i.e., that satisfies
the classical filter equation \ref{classicalorthh}), that in addition 
satisfies $p_0(x)=0$ for $x\in\pm(\frac 17-\epsilon,\frac 3{14}+\epsilon)\cup(\frac 37-\epsilon,\frac 47+\epsilon).$ 
Note that by \ref{classicalorthh}, we then have $p_0(x)=\sqrt 2$ for
$x\in\pm(\frac 27-\epsilon,\frac 5{14}+\epsilon)\cup(-\frac 1{14}-\epsilon,\frac 1{14}+\epsilon)$. 
Let $p_1$ be the standard choice of 
$C^\infty$ high-pass filter associated to $p_0,$ given by $p_1(x)=e^{2\pi ix}\overline{p_0(x+\frac 12)}$.     
Then it is easily checked that the following functions satisfy \ref{orthh}, \ref{orthg} and \ref{orthgh},
and thus are generalized low- and high-pass filters by our definitions in Section 2.  

$$h_{1,1}(x)=\begin{cases} p_0(x)&x\in [-\frac 2{7},\frac 27)\cr
0&\text{\rm \ otherwise}\end{cases} 
$$

$$h_{1,2}(x)=\begin{cases} p_0(x+\frac 12)&x\in [-\frac 1{7},\frac 17)\cr
0&\text{\rm \ otherwise}\end{cases} 
$$

$$h_{2,1}(x)=\begin{cases} \sqrt 2 &x\in \pm [\frac 37,\frac 12)\cr
0&\text{\rm \ otherwise}\end{cases} 
$$

$$h_{2,2}(x)=0$$

$$g_1(x)=\begin{cases}
p_1(x)&x\in[-\frac 2{7},\frac 2{7})\cr
0&\text{\rm \ otherwise}\end{cases} 
$$

$$g_2(x)=\begin{cases} p_1(x+\frac 12)&x\in [\frac {-1}7,\frac 1{7})\cr
0&\text{\rm \ otherwise}\end{cases} 
$$

Now we check that the resulting wavelet does in fact have a $C^{\infty}$ transform.  We first fix an $x$ and show that $\widehat{\phi}_1$ and $\widehat{\phi}_2$
are $C^{\infty}$ in a neighborhood of $x$.  Recall that these functions form the 
first column of the infinite product matrix
$\prod_{j=1}^\infty\frac 1{\sqrt 2} \left(\begin{matrix} h_{1,1}(\frac x{2^j})&h_{1,2}(\frac x{2^j})\\
h_{2,1}(\frac x{2^j})&0\end{matrix} \right)$.  
Because $h_{2,1}$ is 0 on $(-\frac 37, \frac 37)$, all but a finite
number of the matrix factors are upper triangular.  Thus, each of the entries in the first column of the infinite product contains
only a finite number of terms.  Since $h_{1,1}=\sqrt{2}$ in a neighborhood of $0$, each of the terms has only a finite number of factors not equal to 1.  
Thus, it will suffice to show that each of 
the factors in each of the terms is $C^{\infty}$.

By construction, we have that the $h_{i,j}$ are all
$C^\infty$ everywhere except for $h_{1,1}$ at $n\pm\frac 27$, $h_{1,2}$ at $n\pm\frac 17$, and
$h_{2,1}$ at $n\pm\frac 37$ for $n\in\Bbb Z$. We will show that whenever one of these discontinuities occurs as a factor
in the infinite product, it is cancelled by a following factor that is $0$ in a neighborhood of the point of discontinuity.  
Note that the product of two adjacent factors in the 
infinite matrix product has the form
$$\left(\begin{matrix} h_{1,1}(y)h_{1,1}(\frac y2)+h_{1,2}(y)h_{2,1}(\frac y2)&h_{1,1}(y)h_{1,2}(\frac y2)\\
h_{2,1}(y)h_{1,1}(\frac y2)&h_{2,1}(y)h_{1,2}(\frac y2)\end{matrix} \right).$$
Thus any term in the infinite product that contains a factor of $h_{2,1}(n\pm\frac 37)$ must also contain
a factor of one of the forms $h_{1,1}(n\pm\frac 2{7})$, $h_{1,1}(n\pm\frac 3{14})$, $h_{1,2}(n\pm\frac 3{14})$ $h_{1,2}(n\pm\frac 2{7})$.
The first three possibilities are $0$ in a neighborhood of the point in question, so if we have a discontinuous factor of $h_{2,1}$, it is
either cancelled out by a $0$ factor, or we also have a factor of $h_{1,1}(n\pm\frac 27)$ with a smaller $n$.  Similarly,
any term in the infinite product that contains a factor of $h_{1,1}(n\pm\frac 27)$ must also contain
a factor of one of the forms $h_{1,1}(n\pm\frac 5{14})$,  $h_{1,1}(n\pm\frac 1{7})$, 
$h_{1,2}(n\pm\frac 5{14})$, or $h_{1,2}(n\pm\frac 1{7})$. The first three of these possibilities are 0 in a neighborhood of the points in question, so any discontinuous factor of $h_{1,1}$ is either cancelled out by a $0$ factor, or is followed by a factor of $h_{1,2}(n\pm\frac 17)$ with an equal or smaller
$n$.   
Finally, any term in the infinite product that contains a factor
of $h_{1,2}(n\pm\frac 17)$ must also contain
a factor of either the form $h_{2,1}(n\pm\frac 1{14})$ or the form $h_{2,1}(n\pm\frac 3{7})$. The first of these possibilities is $0$ in again 
$0$ in a neighborhood of the point in question; the second possibility throws us back into the first type of discontinuity we considered above, but with a 
smaller $n$.  We can repeat the above sequence of arguments if necessary, noting that each succeeding factor is evaluated at a point half the distance from
the origin as its predecessor, so that the chain above must eventually end with a factor of $0$.   

This argument shows that $\widehat{\phi}_1$ and $\widehat{\phi}_2$
are $C^{\infty}$.  To see that $\widehat{\psi}$ 
is $C^{\infty}$ as well, it suffices to note that $g_1$ and $g_2$ are.

\end{example}

\begin{remark}
With some more work it is possible to slightly improve the arguments above and to show that in fact the Fourier transform vanishes rapidly at infinity.
\end{remark}


\section{A Generalized Loop Groupoid Action on the Bundle of Generalized Filter Systems}

Let $\{h_{i,j}\}_{1\leq\;i,j\;\leq c}$ and 
$\{ g_{k,j}\}_{1\leq\;k\;\leq \tilde{c},\;1\leq\;j\;\leq c}$ be generalized low-pass and high-pass filter functions defined as in Section 2.  
Since $\oplus_{i=1}^c\;L^2(S_i)\;\cong\;L^2(\bigsqcup_{i=1}^c S_i),$ we can suppress the second index of the filter functions and view generalized filter functions as a vector ($c+\tilde{c}$-tuple) of functions: 
$$(h_1,h_2,\;\cdots\;h_c,g_1,g_2,\cdots,g_{\tilde{c}})\;\in\;[L^2(\bigsqcup_{i=1}^c S_i)]^{c+\tilde{c}}.$$
Further, we note that for any fixed $\omega\in \mathbb T^d$, all the information in the output of the vector of functions 
\newline $(h_1,h_2,\;\cdots\;h_c,g_1,g_2,\cdots,g_{\tilde{c}})$ is actually in $\mathbb C^{m(\alpha(\omega))+\tilde{m}(\alpha(\omega))},$ where as in Section 2, $\alpha$ is the endomorphism of $\mathbb T^d$ defined by the matrix $B=A^t,$ since by the orthogonality relations, $h_i(\omega)=0$ if $i>m(\alpha(\omega))$ and  $g_k(x)=0$ if $k>\tilde{m}(\alpha(\omega))$. 

We want to generalize the discussion given by Bratteli and Jorgensen in \cite{BJ1}, and construct a loop groupoid which acts on the class of filter systems with bounded multiplicity functions corresponding to a fixed dilation matrix $A$ acting on $\mathbb R^n.$ We will also impose the condition the filter systems satisfy some ``initial conditions'' that are in some sense canonical.

We first remind readers of the notion of a vector bundle over a topological space $X;$ more  details can be found in \cite{FD}.  

\begin{definition} 
\label{bundledef}
Let $X$ be a topological space.  A (finite dimensional) vector bundle over the space $X,$ denoted by $(E,\;p,\;X),$ is a topological space $E,$ together with a continuous open surjection $p:\; E\;\rightarrow\;X,$ and operations and norms making each fiber $E_x\;=\;p^{-1}(X)$ into a (finite dimensional) vector space, which in addition satisfies the following conditions:
\renewcommand{\labelenumi}{(\roman{enumi})}
\begin{enumerate}
\item $y\mapsto\;\|y\|$ is continuous from $E$ to $\mathbb R,$
\item The operation $+$ is continuous as a function from $\{(y,z)\in E\times E:\;p(y)=p(z)\}$ to $E.$
\item For each $\lambda\;\in\;\mathbb C,$ the map $y\mapsto\;\lambda\cdot y$ is continuous from $E$ to $E.$
\item If $x\in X$ and $\{y_i\}$ is any net of elements of $E$ such that $\|y_i\|\rightarrow 0$ and 
$p(y_i)\rightarrow x$ in $X,$ then $y_i\rightarrow \vec{0}\in E_x$ in $E.$ 
\end{enumerate}
A Borel map $s:X\;\rightarrow\; E$ is called a Borel cross-section if $p\circ s(x)\;=\;x,\forall x\in X.$
\end{definition}
We review the notion of essentially bounded multiplicity functions $m$ associated to a dilation matrix $A$ that can give rise to GMRA's, as described in \cite{BCM} and \cite{BM}.  We first recall that $m$ must satisfy the consistency inequality 
\begin{equation}
\label{consisteq1}
m(\omega)\leq \sum_{l=0}^{N-1}m(\omega_l),
\end{equation}
where we recall from Section 2 that $\{\omega_l:\;0\leq\; l \;\leq N-1\}$ is the set of $N$ distinct preimages of $\omega\in\mathbb T^d$ under the endomorphism $\alpha$ given by $\alpha(\omega)=B(\omega),$ for $\omega\in\mathbb T^d=\mathbb R^d/\mathbb Z^d$ and $B= A^t.$ This inequality leads to the consistency equality Equation \ref{consistency} discussed in Section 2.
Moreover, we can enumerate the $\omega_l$ as follows.  Enumerate a set of coset representatives $\{\xi_0=\vec{0}, \xi_1,\cdots \xi_{N-1}\}$ for $B^{-1}(\mathbb Z^d)/\mathbb Z^d$ including $\vec{0}\in \;B^{-1}(\mathbb Z^d).$  For each $l,\;0\leq l\leq N-1\}$ let $\zeta_l$ be the image of $\xi_l$ under the natural projection from $\mathbb R^d$ onto $\mathbb T^d=\mathbb R^d/\mathbb Z^d.$ 
Note that the $\zeta_l$ are distinct elements of $\mathbb T^d$ and that $\alpha(\zeta_l)=\vec{0},\;0\leq l\leq N-1.$
Find a Borel cross-section $\sigma:\mathbb T^n= \mathbb R^d/\mathbb Z^d\rightarrow \mathbb T^n= \mathbb R^d/B^{-1}(\mathbb Z^d)\cong [\mathbb R^d/\mathbb Z^d]/[B^{-1}(\mathbb Z^d)/\mathbb Z^d],$ with $\sigma(\vec{0})=\vec{0}.$ This Borel map satisfies $\alpha\circ\sigma(\omega)=\omega$ for all $\omega\in\mathbb T^d.$ Then define 
$$\omega_l=\sigma(\omega)\zeta_l,\;0\leq l\leq N-1;$$
one easily verifies that the $\{\omega_l:\;0\leq\; l \;\leq N-1\}$ are the $N$ distinct preimages of $\omega\in\mathbb T^d$ under $\alpha.$

Let $\Delta= \cup_{k=0}^{\infty}B^k(S_1+\mathbb Z^d).$
Recall from \cite{BM} and \cite{BRS} that in order for $m$ to be the multiplicity function for a GMRA, $\Delta$ must satisfy 
\begin{equation}
\label{consisteq2}
\sum_{n\in\mathbb Z^d}\chi_{\Delta}(\omega +n)\geq m(\omega).
\end{equation}
Finally, $\Delta$ should satisfy 
\begin{equation}
\label{consisteq3}
\cup_{p\in\mathbb Z}B^p(\Delta)=\mathbb R^n.
\end{equation}
If conditions \ref{consisteq1}, \ref{consisteq2}, \ref{consisteq3} are satisfied we say that $m$ is an (essentially) bounded multiplicity function with respect to the dilation matrix $A.$
Given such an $m,$ we construct the conjugate multiplicity function $\widetilde{m}$ by defining 
\begin{equation}
\label{conjmult}
\widetilde{m}(\omega)=\sum_{l=0}^{N-1}m(\omega_l)-m(\omega), \omega\in\mathbb T^d;
\end{equation}
by definition, $m$ and $\widetilde{m}$ will satisfy Equation \ref{consistency}.

Given an (essentially) bounded multiplicity function $m$ on $\mathbb T^d,$ let $c\;=\;\text{ess}\;\text{sup}\;m,$ and $\tilde{c}\;=\;\text{ess}\;\text{sup}\;\tilde{m}.$ We recall from Baggett, Courter and Merrill that to such an $m$ we can explicitly construct a canonical system of low-pass filters, or ``generalized conjugate mirror filters'' $\{h^{\mathcal C}_{i,j}:\;1\leq i,j \leq \; c\},$ by using the method of Theorem 3.6 of [BCM]. Given this system $\{h^{\mathcal C}_{i,j}:\;1\leq i,j \; c\},$ Theorem 2.5 item (1) of [BCM] gives us a way of explicitly constructing an associated system of high-pass filters, or ``complementary conjugate mirror filters'', $\{g^{\mathcal C}_{k,j}: 1\leq k \leq \tilde{c},\; 1\leq j \leq c\}$.  We call this family of filters $\{h^{\mathcal C}_{i,j}:\;1\leq i,j\; c\}\cup \{g^{\mathcal C}_{k,j}: 1\leq k \leq \tilde{c},\; 1\leq j \leq c\}$ the {\bf canonical filter system} associated to the multiplicity function $m.$ 

Now let $T_j\;=\{\omega\in\mathbb T^d: m(\alpha(\omega))+\tilde{m}(\alpha(\omega))=j\},\;0\;\leq\;j\;\leq\;c+\tilde{c}.$  Set $T_{i,j}\;=\;S_i\cap T_j,\;0\;\leq\;j\;\leq\;c+\tilde{c};$ then each $T_{i,j}$ is Borel and $S_i\;=\;\bigsqcup_{j=0}^{c+\tilde{c}} T_{i,j}.$
\begin{definition}
\label{Msystemdef}
Fix a bounded multiplicity function $m$ associated  to the dilation matrix $A.$ Let $E_m$ be the Borel space given by 
$$E_m\;=\;\bigsqcup_{i=1}^c \bigsqcup_{j=0}^{c+\tilde{c}} [T_{i,j}\times\mathbb C^j].$$
Let $(E_m,\;p,\;\bigsqcup_{i=1}^c S_i)$ be the Borel vector bundle where the map $p:\;E_m\;\rightarrow\;\:\bigsqcup_{i=1}^c S_i$ is defined by $p(x,\vec{v})=x,\;(x,\vec{v})\;\in\;T_{i,j}\times\mathbb C^j.$ By definition, $(E_m,\;p,\;\bigsqcup_{i=1}^c S_i)$ is a vector bundle over $\bigsqcup_{i=1}^c S_i$ whose fiber over $\omega \in S_i$ is a complex vector space of dimension $m(\alpha(\omega))+\tilde{m}(\alpha(\omega)).$ An {\bf $M$-system associated to the multiplicity function $m$} is a
Borel cross-section $M:\;\bigsqcup_{j=1}^c S_j\;\rightarrow\;E_m$ of this bundle whose values, $(M_1(\omega), M_2(\omega),\;\cdots M_{m(\alpha(\omega))+\tilde{m}(\alpha(\omega))})$, are the output of a vector of generalized low- and high-pass filters. That is, for fixed $\omega\in\bigsqcup_{j=1}^c S_j,$ the first
$m(\alpha(\omega))$ components of the vector $M(\omega)$ correspond exactly to the generalized low-pass filter function values $h_1(\omega), h_2(\omega),\cdots, h_{m((\alpha(\omega))}(\omega),$ and the last $\tilde{m}(\alpha(\omega))$ components of $M(\omega)$ correspond to the generalized high-pass filter function values $g_1(\omega), g_2(\omega),\cdots, g_{\tilde{m}(\alpha(\omega))}(\omega),$ where the filters involved satisfy the low-pass and Lipschitz conditions defined in Section 2, and where in addition the filters $\{h_{i}\}\cup \{g_k\}$ satisfy the initial conditions $h_i(\zeta_l)\;=\;h^{\mathcal C}_i(\zeta_l),\; 1\leq i\leq c,\;0\leq l \leq N-1,$ and $g_k(\zeta_l)\;=\;g^{\mathcal C}_k(\zeta_l), \; 1\leq k\leq \tilde{c},\; 0\leq l \leq N-1,$ where the $h^{\mathcal C}_i$ and the $g^{\mathcal C}_k$ are the canonical filters associated to $m$ in the previous paragraph.  Denote by ${\mathcal M}_m$ the set of $M$-systems associated to the bounded multiplicity function $m.$ We remark that ${\mathcal M}_m$ can be given the structure of a topological space if its elements are viewed as elements of the Hilbert space $\oplus_{i=0}^c\oplus_{j=0}^{c+\tilde{c}}L^2(T_{i,j})\otimes\mathbb C^j.$
\end{definition}

\begin{remark}
The definition above is essentially the same as Definition 4.2 of 
\cite{BJMP} except for the additional assumption that the $M$-system have the appropriate canonical values at the preimages of $\vec{0}$ under the endomorphism $\alpha,$ which was missing from \cite{BJMP}.
Note all the information about the generalized filters $\{h_{i,j}\}$ and $\{g_{k,j}\}$ is encoded in the $M$-system.  In particular, for any fixed multiplicity function $m$, such that both $m$ and $\tilde{m}$ are constant in a neighborhood of the origin, we have a one-to-one correspondence between $M$-systems and   
collections of generalized filter functions as defined in Section 2 which satisfy the canonical initial conditions.  These generalized filters will in turn give rise to an orthonormal frame wavelet family and its associated GMRA $\{V_j\}$ by Theorem \ref{waveletGMRAth} from Section 3.  However, as shown by one of the examples given in Section 2, the multiplicity function $m'$ for the GMRA $\{V_j\}$ {\bf need not} be equal to $m.$  So, there is no correspondence in general, between the class of $M$-systems associated to a given multiplicity function $m$ and the class of GMRA's associated to $m.$
\end{remark}

To develop the loop group action on these $M$-systems, we first define an endomorphism $\Pi_{\alpha}:\;\bigsqcup_{i=1}^c S_i\;\rightarrow\;\mathbb T^d$ by $\Pi_{\alpha}(\omega)\;=\;\alpha(\omega).$  Each $\omega\in\mathbb T^d$ has $\sum_{l=0}^{N-1}m(\omega_l)\;=\;m(\omega)\;+\;\tilde{m}(\omega)$ preimages in $\bigsqcup_{i=1}^c S_i,$ where the $\{\omega_l: 0\leq l\leq N-1\}$ are the $N$ distinct preimages of $\omega$ in $\mathbb T^d$ under the endomorphism $\alpha,$ parametrized as discussed earlier.
For convenience of notation, we label these preimage maps $r_{(l,j)},$ where $r_{(l,j)}(\omega)\;=\;\omega_l\;\in\;S_j\;\subseteq\;\bigsqcup_{i=1}^c S_i$ for 
$1\leq\;j\;\leq m(\omega_l).$   (Note that this range on $j$, as $l$ varies from 0 to $N-1$, gives all the preimages, since if $j> m(\omega_l),$ by definition $\omega_l$ is not an element of $S_j.$) For each fixed $\omega$, we give the pairs $(l,j)$ the lexicographical order, so that $(l_1,j_1) \leq (l_2,j_2)$ if $l_1 < l_2$ or if $l_1 = l_2$ and $j_1 \leq j_2.$  We thus implicitly define a 1-1 map $\lambda_{\omega}$ taking the pairs $(l,j)$ onto the integers from 1 to $m(\omega)+\tilde{m}(\omega)$. 

We now construct a unitary group bundle $(F,\; q,\;\mathbb T^d)$ as follows.  For each $j\;\in\;\{1,\cdots,c+\tilde{c},\},$ let 
$Z_j\;=\;\{\omega\;\in\;\mathbb T^d:\;m(\omega)+\tilde{m(\omega)}=j\}.$ Let 
$$F_m\;=\;\bigsqcup_{j=0}^{c+\tilde{c}} [Z_j\times U(j,\mathbb C)],$$
where $U(j,\mathbb C)$ is viewed as a the topological group of unitary $j\times j$ matrices:  it inherits its topology from being a subset of $C^{\ast}$-algebra of complex $j\times j$ matrices given the operator norm. 
Define $q:\; F_m\;\rightarrow\;\mathbb T^d$ by $q(\omega,T)\;=\;\omega,$ for $(\omega,T)\;\in\; F_m,$ and note that  $q:\;F_m\;\rightarrow\;\mathbb T^d$ is a continuous open surjection, and the fiber $q^{-1}(\omega)$ of the bundle consists of the group of complex unitary matrices $U(m(\omega)+\tilde{m}(\omega),\mathbb C).$ Borel cross sections to this group bundle consist of Borel maps $L:\mathbb T^d\;\rightarrow\;F_m$ such that $q\circ L(\omega)\;=\;\omega.$
Note also that $(F_m,\;q,\mathbb T^d)$ is a subset of a Borel vector bundle over $\mathbb T^d$ in the sense of Definition \ref{bundledef}, whose fibers over particular values of $\omega\in\mathbb T^d$ are made up of algebras of finite-dimensional matrices of varying dimension.
We denote the set of sections of this bundle by $\Gamma_m(F_m,q).$ Note $\Gamma_m(F_m,q)$ is a group under pointwise operations on $\mathbb T^d,$ where the identity element of the group is given by that section whose value at $\omega$ is equal to $Id_{m(\omega)+\tilde{m}(\omega)}.$ 

Also, it is possible to define a groupoid corresponding to a dilation matrix $A$ as follows.  Let $\Omega_A$ consist of the set of all bounded multiplicity functions corresponding to the dilation matrix $A,$ that is, 
$$\Omega_A\;=\;\{m:\mathbb T^d\rightarrow \mathbb N\cup\{0\}:\;m\;\text{is essentially bounded and satisties}\; \ref{consisteq1}, \ref{consisteq2}, \ref{consisteq3}\},$$
where two multiplicity functions for $A$ are identified if they are equal almost everywhere on $\mathbb T^d.$ There are a variety of topologies we can put on $\Omega_A;$ for the moment the Hilbert space topology of $L^2$ functions on $\mathbb T^d$ will do.  We call $\Omega_A$ {\bf the multiplicity function space associated to} $A.$  Now set 
$${\mathcal L}_A\;=\;\bigsqcup_{m\in\Omega_A}\Gamma_m(F_m,q).$$ 
Then ${\mathcal L}_A$ is a groupoid, with range map $r$ equal to the source map $s$ defined from 
${\mathcal L}_A$ to $\Omega_A$ by $r(\gamma)=s(\gamma)=m$ for $\gamma\in\Gamma_m(F_m,q).$
Hence two elements $\gamma_1$ and $\gamma_2$ of ${\mathcal L}_A$ can be multiplied together if and only if $\gamma_1,\;\gamma_2\;\in\;\Gamma_m(F_m,q)$ for a fixed bounded multiplicity function $m.$

We are now state a theorem about $M$-systems that can be derived in a fairly straightforward fashion from the orthogonality relations.  This theorem is a generalization of Theorem 4.3 of 
\cite{BJMP}, and thus we merely sketch the proof and refer the reader to \cite{BJMP} for details.
\begin{theorem} 
\label{mainthm}
Let  $\Gamma_m(F_m,q)$ be the group of cross sections of the group bundle associated to a fixed bounded multiplicity function $m$ defined above.  Let $M:\;\bigsqcup_{j=1}^c S_j\;\rightarrow\;E$ be an $M$-system associated to $m.$ Then 
 $\omega\;\mapsto\;(L_{i,\lambda_{\omega}(l,j)}(\omega)),$ where 
$$L_{i,\lambda_{\omega}(l,j)}(\omega)\;=\;\sqrt{\frac{1}{N}}M_i(r_{(l,j)}(\omega))$$
is an element of $\Gamma_m(F_m,q).$ 
\end{theorem}
\begin{proof}
As noted above, $1\leq\lambda_{\omega}(l,j)\leq m(\omega)+\tilde{m}(\omega)$, so for each $\omega\in\mathbb T,$ the matrix $(L_{i, \lambda_{\omega}(l,j)}(\omega))$ is a square matrix of the correct dimension. 
 
We shall show that for all $\omega\in \mathbb T^d,$ the rows of $(L_{i,\lambda_{\omega}(l,j)}(\omega))$  are orthonormal, by means of the orthogonality relations for generalized filter functions given in Equation \ref{orthcol} in the statement of Theorem \ref{thmorthocol}. 

Write $L_i$ for the $i$th row of $(L_{i,\lambda_{\omega}(l,j)}(\omega)).$ If $1\;\leq i\leq i'\leq\; m(\omega),$  
\begin{eqnarray}
<L_i,\;L_{i'}>&=&\sum_{\lambda_{\omega}(l,j)=1}^{m(\omega)+\tilde{m}(\omega)}L_{i,\lambda_{\omega}(l,j)}(\omega)\overline{L_{i,\lambda_{\omega}(l,j)}(\omega)} \nonumber \\
&=& \sum_{l=0}^{N-1}\sum_{j=1}^{m(\omega_l)}\sqrt{\frac{1}{N}}M_i(r_{(l,j)}(\omega))\overline{\sqrt{\frac{1}{N}}M_{i'}(r_{(l,j)}(x))} \nonumber \\
&=& \sum_{l=0}^{N-1}\sum_{j=1}^{m(\omega_l)}\frac{1}{N}h_{i,j}(\omega_l)\overline {h_{i',j}(\omega_l)} \nonumber \\
&=& \frac{1}{N}\sum_{j=1}^{c}\sum_{l=0}^{N-1}h_{i,j}(\omega_l)\overline {h_{i',j}(\omega_l)} \nonumber
\end{eqnarray}
(as $h_{i,j}(\omega_l)\;=\;0$ for $j>\mu(\omega_l),$ since $\omega_l\notin\;S_j$ in that case)  
$$\;=\;(\text{by Equation (2.2)})\;\;N\frac{1}{N}\delta_{i,i'}\chi_{S_i}(\omega)\;=\;\delta_{i,i'}$$
(we note that $\chi_{S_i}(\omega)=1$ since we have $i\leq\; m(\omega)$ and for those values of $i,\;\omega\;\in\;S_i$ by definition of $m(\omega)$).

The cases $m(\omega)\;< i\leq i'\leq\; m(\omega)+\tilde{m}(\omega)$ and $1\;\leq i \leq\;m(\omega)<\; i'\leq\;m(\omega)+\tilde{m}(\omega)$  follow from similar arguments using Equations 2.4 and 2.5.  Thus we have that in all cases, the rows of $(L_{i,\lambda_{\omega}(l,j)}(\omega))$ are orthonormal, and we have the desired unitary matrix.
\end{proof}

\begin{remark}
\label{initialcondition}
We note that for $\omega\in\mathbb T^d,$ the $(m(\omega)+\tilde{m}(\omega))\times(m(\omega)+\tilde{m}(\omega))$ matrix $L(\omega)$ is exactly the submatrix $L(\omega)$ of the matrix $K(\omega)$ defined in Theorem 2.5.  Thus in some sense the proof given above is redundant. We have added in the extra detail because we want to exercise care in enumerating the rows and columns of $L(\omega)$ for future use. We also note that the ``initial condition'' on $M$-systems given in Definition \ref{Msystemdef} can be rephrased as follows.  Let $L^{\mathcal C}_{i,\lambda_{\omega}(l,j)}(\omega))$ be the element of    $\Gamma_m(F_m,q)$ corresponding to the ``canonical'' filter system associated to the multiplicity function $m,$ and let $L_{\mathcal C}$ be  be the $[m(\vec{0})+\tilde{m}(\vec{0})]\times [m(\vec{0})+\tilde{m}(\vec{0})]$ unitary matrix $L^{\mathcal C}_{i,\lambda_{\omega}(l,j)}(\vec{0}).$ Then for any $M$ system $M$ associated to $m,$ we have $\sqrt{\frac{1}{N}}M_i(r_{(l,j)}(\vec{0}))\;=\;L_{\mathcal C}.$  This follows from the initial conditions satisfied by any $M$ system outlined in Definition \ref{Msystemdef}. 
\end{remark}

The results of Theorem \ref{mainthm} imply that the columns of  $(L_{i,\lambda_{\omega}(l,j)}(\omega))$ are orthonormal as well, and allow us to obtain the following corollary, which is a generalization of Corollary 4.4 of \cite{BJMP}. 
\begin{corollary} 
\label{cororthcol}
Let $m$ and $\tilde{m}$ be bounded multiplicity and conjugate multiplicity functions associated to a dilation matrix $A$ that are constant in a neighborhood of the origin, with related sequences of sets 
$\{S_i\:|\;1\leq\; i\;\leq c\}$ and $\{\tilde{S}_k\:|\;1\leq\; k\;\leq \tilde{c}\}.$ Suppose $\{h_{i,j}\}_{1\leq\;i,j\;\leq c}$ and 
$\{ g_{k,j}\}_{1\leq\;k\;\leq \tilde{c},\;1\leq\;j\;\leq c}$ are generalized low-pass and high-pass filter functions associated to the multiplicity function $m.$   
Then for all $\omega\;\in\;\mathbb T^d,$ and for all $l,\;l'\;\in\;\{0,1,\cdots,N-1\}$ and $j,\;j'\;\in\;\{1,\cdots,c)\},$ we have 
$$\sum_{i=1}^{c+\tilde{c}}\frac{1}{N}M_i(r_{(l,j)}(\omega))\overline{M_i(r_{(l',j')}(\omega))}$$
$$\;=\;\sum_{i=1}^c \frac{1}{N}h_{i,j}(\omega_l)\overline{h_{i,j'}(\omega_l)}\;+\;
\sum_{k=1}^{\tilde{c}} \frac{1}{N}g_{k,j}(\omega_l)\overline{g_{k,j'}(\omega_l)}\;=\;
\delta_{j,j'}\delta_{l,l'},$$ 
where here the $M_i$ correspond to the generalized filter functions $h_{i,j}$ and $g_{k,j}$ as given in Definition \ref{Msystemdef}. 
\end{corollary}
\begin{proof}
This result follows fairly directly from Theorems 2.5 and \ref{mainthm}, and we omit the details of the proof. 
\end{proof}

We are ready to define the generalized loop group and its associated action on the set of $M$-systems associated to a multiplicity function $m$.
\begin{definition} 
The {\bf loop group} associated to the multiplicity function $m$ 
is defined to be the subgroup $\text{Loop}_m(F_m,q)$ of the group of Borel sections $\Gamma_m(F_m,q)$ whose elements $K$ 
satisfy $L(\vec{0})\;=\;Id_{m(\vec{0})+\tilde{m}(\vec{0})},$ and  $L_{i,\lambda_{\omega}(j,l)}$ are $\text{Lip}_1$ in a neighborhood of the origin.  
\end{definition}

We now prove that the generalized loop group $\text{Loop}_m(F,q)$ acts freely and transitively on the set of $M$-systems associated to the multiplicity function $m$.  The following theorem generalizes Theorem 4.7 of \cite{BJMP}.
\begin{theorem}
\label{loopthm}
There is a free and transitive action of $\text{Loop}_m(F,q)$ on the set of $M$-systems associated to an essentially bounded multiplicity function $m$ such that $m$ is constant in neighborhoods of $0_l,\;0\leq l\leq N-1,$ where the set $\{0_l:\;0\leq l\leq N-1\}$ consists of the preimages of $0$ under the endomorphism $\alpha:\mathbb T^d\rightarrow\; \mathbb T^d$ This action is given by 
$$L\cdot M(\omega)\; =\; L(\Pi_{\alpha}(\omega))[(M_1(\omega),M_2(\omega),\cdots, M_{m(\alpha(\omega))+\tilde{m}(\alpha(\omega))}(\omega))]^t.$$
\end{theorem}
\begin{proof}
We first prove the transitivity.  Suppose we are given two different $M$-systems, labeled $M=(M_i)$ and $\widetilde{M}=(\widetilde{M}_i).$ Define an element $L$ of the group bundle associated to $m$, that is, an element of $\Gamma_m(F,q),$ 
where $L(\omega)$ has dimension $m(\omega)+\tilde {m}(\omega)$, as follows:
\begin{equation}
\label{eqtransitive}
L_{i,i'}(\omega)=\frac{1}{N}\sum_{\lambda_{\omega}(l,j)=1}^{m(\omega)+\tilde{m}(\omega)} \overline {M_{i'}(r_{(l,j)}(\omega))}\widetilde{M_i}(r_{(l,j)}(\omega)).
\end{equation}  
Since the $M$-systems $M$ and $\widetilde{M}$ have the same values at $\vec{0}$ and at the preimages $\zeta_l$ of $\vec{0}$ under the automorphism $\alpha,$ one easily verifies that $L(\vec{0})\;=\;Id_{m(\vec{0})+\tilde{m}(\vec{0})}.$
Also, one sees by inspection that if $L_M$ and $L_{\widetilde{M}}$ are the elements of $\Gamma_m(F_m,q)$ corresponding to the $M$-systems $M$ and $\widetilde{M}$ as in Theorem \ref{mainthm}, then the proposed element $L$ of $\text{Loop}_m(F,q)$ given in Equation \ref{eqtransitive} is exactly $L(\omega)\;=\;L_{\widetilde{M}}(\omega)[L_{M}(\omega)]^{\ast}.$ 

In addition, as shown by Bratteli and Jorgensen in the classical case, we have 
\begin {eqnarray}
[L\cdot M]_i(\omega)&=& \sum_{i'=1}^{m(\alpha(\omega))+\tilde{m}(\alpha(\omega))}L_{i,i'}(\Pi_{\alpha} (\omega))M_{i'}(\omega) \nonumber \\
&=& \sum_{i'=1}^{m(\alpha(\omega))+\tilde{m}(\alpha(\omega))}\left(\frac{1}{N}\sum_{\lambda_{\alpha(\omega)}(l,j)=1}^{m(\alpha(\omega))+\tilde{m}(\alpha(\omega))} \overline {M_{i'}(r_{(l,j)}(\alpha(\omega)))}\widetilde{M_i}(r_{(l,j)}(\alpha(\omega)))\right)M_{i'}(\omega) \nonumber \\
&=& \sum_{\lambda_{\alpha(\omega)}(l,j)=1}^{m(\alpha(\omega))+\tilde{m}(\alpha(\omega))}\widetilde{M_i}(r_{(l,j)}(\alpha(\omega)))\left(\sum_{i'=1}^{m(\alpha(\omega))+\tilde{m}(\alpha(\omega))}\overline {M_{i'}(r_{(l,j)}(\alpha(\omega)))}M_{i'}(\omega)\right) \nonumber \\
&=& \widetilde{M_i}(\omega), \nonumber
\end {eqnarray}
where the last equality follows from the orthogonality of the columns of $M$ as established in Corollary \ref{orthcol}, so that the inside sum is 0 except for the single values of $l$ and $j$ where $r_{l,j}(\alpha(\omega))=\omega$. 

To prove that the action is free, suppose $M=(M_i)$ is a $M$-system associated to $m$ and $L\;\in\;\text{Loop}_m(F,q)$ satisfies 
$$L(\Pi_{\alpha}(\omega))[(M_1(\omega),M_2(\omega),\cdots, M_{m(\alpha(\omega))+\tilde{m}(\alpha(\omega))}(\omega))]^t\;=$$
$$[(M_1(\omega),M_2(\omega),\cdots, M_{m(\alpha(\omega))+\tilde{m}(\alpha(\omega))}(\omega))]^t.$$
Using Theorem \ref{mainthm}, for each $\omega\in\mathbb T^d$ we define a $(m(\alpha(\omega))+\tilde{m}(\alpha(\omega)))\times (m(\alpha(\omega))+\tilde{m}(\alpha(\omega)))$ unitary matrix ${\mathcal M}$ by 
$${\mathcal M}_{i,\lambda_{\alpha(\omega)}(l,j)}(\omega)\;=\;\sqrt{\frac{1}{N}}M_i (r_{(l,j)}(\alpha(\omega))),$$
$$1\;\leq\; i\;\leq\; m(\alpha(\omega))+\tilde{m}(\alpha(\omega)),\;0\leq\;l\;\leq\;N-1,\;1\leq\;j\;\leq m(\alpha(\omega)_l)).$$  
We then see that $L(\alpha(\omega)){\mathcal M}(\omega)\;=\;{\mathcal M}(\omega)$ for all $\omega\in\;\mathbb T^d.$ 
By unitarity of ${\mathcal M}(\omega),$ this shows that  that $L(\alpha(\omega))$ is the $(m(\alpha(\omega))+\tilde{m}(\alpha(\omega)))\times (m(\alpha(\omega))+\tilde{m}(\alpha(\omega)))$ identity matrix for all $\omega\in\;\mathbb T^d,$ which implies that $L$ is the identity element of $\text{Loop}_m(F,q),$ as desired. 
\end{proof}

Similarly, we can define the generalized loop groupoid assigned to a dilation matrix $A.$
\begin{definition}
Let $A$ be a $d\times d$ integer dilation matrix. Let $\Omega_A$ be the multiplicity function space associated to $A$.  The {\bf loop groupoid} associated to $A$ consists of the set 
$${\mathcal L}oop_A\;=\;\bigsqcup_{m\in\Omega_A}\text{Loop}_m(F_m,q),$$
where $r=s:{\mathcal L}oop_A\rightarrow \Omega_A$ is defined by $r(L)=s(L)=m$ for $L \in\text{Loop}_m(F_m,q).$
\end{definition}

We can also define a fiber bundle of $M$-systems over $\Omega_A$ as follows:
\begin{definition} Let $\Omega_A$ be the multiplicity function space associated to $A$.
Define ${\mathbb M}_A\;=\cup_{m\in \Omega_A}\{m\}\times{\mathbb M}_m,$ where ${\mathbb M}_m$ consists of the set of $M$-systems associated to a fixed multiplicity function $m,$ where two $M$-systems are identified if they are equal almost everywhere.  Give elements of ${\mathbb M}_m$ the Hilbert space topology mentioned in Definition \ref{Msystemdef}  The topology on ${\mathbb M}_A$ is obtained from viewing an element $(m,{\mathcal M})\in{\mathbb M}_A$ as a tuple $(m,{\mathcal M})$ of functions defined on the Cartesian product space $\mathbb T^d\times \bigsqcup_{i=1}^{\infty}[\mathbb T^d]_i,$ taking on values in $[\mathbb N\cup\{0\}]\oplus l_2(\mathbb N),$ where 
$$(m,{\mathcal M})(\omega,z)=(m(\omega),{\mathcal M}_1(z),{\mathcal M}_2(z),\cdots {\mathcal M}_{m(\alpha(z))+\tilde{m}(\alpha(z))}(z),0,0,0,\cdots),$$ 
and where ${\mathcal M}_j(z)=M_j(z)$ if $z\in\bigsqcup_{i=1}^c S_i$ and $1\leq j\leq m(\alpha(z))+\tilde{m}(\alpha(z)),$ and 
${\mathcal M}_j(z)=0$ if $z\notin\bigsqcup_{i=1}^c S_i,$ for $\{S_i\}$ the standard subsets of $\mathbb T^d$ associated to the multiplicity function $m.$  Then ${\mathbb M}_A$ is a topological space, if we view elements of ${\mathbb M}_A$  as elements in the Hilbert space $L^2(\mathbb T^d)\oplus L^2(\bigsqcup_{i=1}^{\infty}[\mathbb T^d]_i)\otimes[ l_2(\mathbb N)]$.
Define a map $\Pi:{\mathbb M}_A\rightarrow \Omega_A$ by 
$\Pi((m,{\mathcal M}))=m.$ Then $\Pi$ is a continuous surjection, since it is a restriction of the Hilbert space projection from $L^2(\mathbb T^d)\oplus L^2(\bigsqcup_{i=1}^{\infty}[\mathbb T^d]_i)\otimes[l_2(\mathbb N)]$ onto $L^2(\mathbb T^d)\oplus [\{\vec{0}\}]\cong L^2(\mathbb T^d),$ and $\Pi^{-1}(m)= {\mathbb M}_m,$ so that ${\mathbb M}_A$ is a fiber bundle over $\Omega_A,$ called the {\bf fiber bundle of $M$-systems associated to $A$}.   
\end{definition} 

Using these definitions, we obtain the following corollary to Theorem \ref{loopthm}.
\begin{corollary}
Let $A$ be a $d\times d$ integer dilation matrix, and let ${\mathcal L}oop_A$ and ${\mathbb M}_A$ the loop groupoid and fiber bundle of $M$-systems associated to $A,$ respectively.  Then there is a groupoid action of ${\mathcal L}oop_A$ on ${\mathbb M}_A,$ where $L \in{\mathcal L}oop_A$ is allowed to act on ${\mathcal M}\in{\mathbb M}_A$ if and only if $r(L)=\Pi({\mathcal M})=m.$ This action is fiberwise transitive.
\end{corollary}
\begin{proof}
This is just a restatement of part of Theorem \ref{loopthm}.
\end{proof}

We now modify Example 4.5 of \cite{BJMP}, in order to show how we can use the loop group action to transform the canonical filter functions for the Journ\'e wavelet into the filter functions discussed in Example \ref{ExJourne}.  The construction of the canonical filter functions for the Journ\'e wavelet was first done in the thesis of J. Courter \cite{Cou}. 
\begin{example}
The Journ\'e wavelet in the frequency domain is 
the characteristic function of the set 
$$[-\frac{16}{7},-2)\cup [-\frac{1}{2},-\frac{2}{7})\cup [\frac{2}{7},\frac
{1}{2}]\cup [2,\frac{16}{7}).$$
Here the multiplicity function $m$ takes on the values $0,1,$ and $2,$ and $\tilde{m}(x)\equiv 1,$ since the Journ\'e wavelet is a single orthonormal wavelet.
If we identify $\mathbb T$ with  $[-\frac{1}{2},\frac{1}{2}),$ we can write  
$S_1=[-\frac {1}{2},-\frac{3}{7})\cup [-\frac {2}{7},\frac {2}{7})\cup 
[\frac {3}{7},\frac{1}{2}),\;S_2=[-\frac {1}{7},\frac{1}{7}),$ and 
$\widetilde{S_1}=[-\frac {1}{2},\frac{1}{2}].$
The canonical generalized filter functions then are: 
$$h^{\mathcal C}_{1,1}(x)=\sqrt{2}\chi_{[-\frac{2}{7},-\frac{1}{4})\cup [-\frac{1}{7},\frac{1}{7})\cup 
[\frac{1}{4},\frac{2}{7})}(x),$$
$$h^{\mathcal C}_{2,1}=\sqrt{2}\chi_{[-\frac{1}{2},-\frac{3}{7})\cup [\frac{3}{7},\frac{1}{2})}(x),$$
$$g^{\mathcal C}_1(x)=\sqrt{2}\chi_{[-\frac{1}{4},-\frac{1}{7})\cup [\frac{1}{7},\frac{1}{4})}(x);$$
$$h^{\mathcal C}_{1,2}(x)=0,$$
$$h^{\mathcal C}_{2,2}(x)=0,$$
$$g^{\mathcal C}_2(x)=\sqrt{2}\chi_{[-\frac{1}{7},\frac{1}{7})}(x).$$
One calculates that the ``initial condition'' matrix $L_{\mathcal C}$ discussed in Remark \ref{initialcondition} corresponding to this canonical filter system is the $3\times 3$ matrix 
$$L_{\mathcal C}\;=\;\left(\begin{array}{ccc}
1\;&0\;&0\\
0\;&0\;&1\\
0\;&1\;&0\end{array}\right).$$
\vskip.1in
\noindent  We now construct an element of the loop group $L_p$ such that if $M_J$ is the $M$-system corresponding to the above output of generalized low- and high-pass filters, the $M$-system $L_p\cdot M_J$ has as its output the filter functions corresponding to those given in the example in Section 3.
\noindent Consider the decomposition of the disjoint union $S_1\bigsqcup S_2 $ (identified with a subset of  $[-\frac{1}{2},\frac{1}{2})\bigsqcup [-\frac{1}{2},\frac{1}{2})$) 
given by
\vskip.1in

$T_{1,1}= \pm  [\frac{1}{7},\frac{3}{14})$  (Here $m(x)\geq 1,$ $m(2x)+\tilde{m}(2x)=1.)$

$T_{1,2}= \pm [\frac{1}{14},\frac{1}{7})\cup \pm [\frac{3}{14},\frac{2}{7})  $ (Here $m(x)\geq 1,$ $m(2x)+\tilde{m}(2x)=2.)$

$T_{1,3}= [-\frac{1}{14},\frac{1}{14})\cup \pm [\frac{3}{7},\frac{1}{2})$  (Here $m(x)\geq 1,$ $m(2x)+\tilde{m}(2x)=3.)$

$T_{2,2}=  \pm [\frac{1}{14},\frac{1}{7})$  (Here $m(x)=2,$  $m(2x)+\tilde{m}(2x)=2.)$

$T_{2,3}=[-\frac{1}{14},\frac{1}{14})$ (Here $m(x)=2,$ $m(2x)+\tilde{m}(2x)=3.)$

Then as in the general case, $S_j\;=\;\bigsqcup_{k=0}^{3} T_{j,k},\;j=1,2.$

\vskip.1in
\noindent We now describe the $M$-system $M_J:\bigsqcup_{i=1}^2 S_i\;\rightarrow E_m$ associated to the filter functions above.  We describe $M_J$ as a cross-section separately on both of the (disjoint) sets $S_1$ and $S_2.$  

On $S_1,\;M_J$ is given by: 
$$[M_J(x)]=\left\{\begin{array}{rrrrr}
{[\sqrt{2}],}&\mbox{if}\ \;x\in T_{1,1},\\
{[\sqrt{2},0], }&\mbox{if}\ \;x\in \pm [\frac{1}{14},\frac{1}{7})\cup \pm [\frac{1}{4}, \frac{2}{7}) \subseteq T_{1,2},\\
{[0,\sqrt{2}], }&\mbox{if}\ \;x\in \pm [\frac{3}{14},\frac{1}{4}) \subseteq T_{1,2},\\
{[0,\sqrt{2},0],}&\mbox{if}\ \;x\in T_{1,3}\backslash T_{2,3},\\
{[\sqrt{2},0,0],}&\mbox{if}\;x\in [T_{1,3}\cap T_{2,3}].
\end{array}\right.$$

\vskip.1in
On $S_2,\;M_J$ is given by: 
$$[M_J(x)]=\left\{\begin{array}{rr}
{[0,\sqrt{2}],}&\mbox{if}\;x\in T_{2,2},\\
{[0,0,\sqrt{2}],}&\mbox{if}\;x\in T_{2,3}.
\end{array}\right.$$

\vskip.1in

 Consider the decomposition of the circle $\mathbb T$ (identified with $[-\frac{1}{2},\frac{1}{2})$) 
given by
\vskip.1in

$P_1=[-\frac{1}{7},\frac{1}{7})$ (Here $m(x)=2, m(\frac{x}{2})=2, m(\frac{x+1}{2})=1.)$

$P_2=\pm [\frac{1}{7},\frac{2}{7})$ (Here $m(x)=1,$ $m(\frac{x}{2})=2,$ $m(\frac{x+1}{2})=0.)$

$P_3=\pm [\frac{2}{7},\frac{3}{7})$, (Here $m(x)=0$, $m(\frac{x}{2})=1,\; m(\frac {x+1}2)=0.$)

$P_4=\pm [\frac{3}{7},\frac{1}{2})$  (Here $m(x)=1,$ $m(\frac{x}{2})=1,\; m(\frac {x+1}2)=1.$)

We note the following, which will be useful in our calculations: $2T_{1,1} = P_3,
\newline 2(\pm [\frac{3}{14},\frac{1}{4}))\;= 2(\pm [\frac{1}{4},\frac {2}{7}))= P_4,\;2T_{1,3}=P_1,\; 2 (\pm [\frac{1}{14},\frac{1}{7})) = 2T_{2,2}=P_2,$ and $2T_{2,3} = P_1.$

\vskip.1in
Now as in Example \ref{ExMult1}, let  $p_0$ be any classical (MRA)  low-pass filter for dilation by $2,$ (\i.e., one that satisfies
the classical filter equation \ref{classicalorthh}), that also 
satisfies $p_0(x)=0$ for $x\in\pm(\frac 17-\epsilon,\frac 3{14}+\epsilon)\cup(\frac 37-\epsilon,\frac 47+\epsilon).$  Note that by Equation \ref{classicalorthh}, we then have $p_0(x)=\sqrt 2$ for
$x\in\pm(\frac 27-\epsilon,\frac 5{14}+\epsilon)\cup(-\frac 1{14}-\epsilon,\frac 1{14}+\epsilon)$.
Let $\widetilde{p_0}=\frac{p_0}{\sqrt{2}}$ be the normalization of $p_0,$ so that  $\widetilde{p_0}(0)=1.$ Let $p_1$ be the standard choice of associated high-pass filter, so that $p_1(x)=e^{2\pi ix}\overline {p_0(x+\frac{1}{2})},$ and let $\widetilde{p_1}=\frac{p_1}{\sqrt{2}}$ be the normalization of $p_1.$ Again, by Equation \ref{classicalorthh},  $\widetilde{p_1}$ must have modulus $1$ for $x\in\pm(\frac 17-\epsilon,\frac 3{14}+\epsilon)\cup(\frac 37-\epsilon,\frac 47+\epsilon).$ 

\vskip.1in
The associated element of the loop group $L_{p}$ that we want to choose is:
$$L_{p}(x)=\left\{\begin{array}{ll}
{\left(\begin{array}{ccc}
\widetilde{p_0}(\frac{x}{2})\;&\;0&\widetilde{p_0}(\frac{x+1}{2})\\
0\;&1\;&0\\
\widetilde{p_1}(\frac{x}{2})\;&\; 0&\widetilde{p_1}(\frac{x+1}{2})\end{array}\right),}&\mbox{if}\ \;x \in P_1,\\
{\left(\begin{array}{cc}
\widetilde{p_0}(\frac{x}{2})\;&\widetilde{p_0}(\frac{x+1}{2})\\
\widetilde{p_1}(\frac{x}{2})\;&\widetilde{p_1}(\frac{x+1}{2})\end{array}\right),}&\mbox{if}\ \;x\in P_2,\\
\;\;\;\;\;\;\;\;\widetilde{p_1}(\frac{x}{2}),&\mbox{if}\ \;x\in P_3,\\
{\left(\begin{array}{cc}
\widetilde{p_0}(\frac{x+1}{2})\;&\widetilde{p_0}(\frac{x}{2})\\
\widetilde{p_1}(\frac{x+1}{2})\;&\widetilde{p_1}(\frac{x}{2})\end{array}\right),}&\mbox{if}\ \;x\in P_4.\\
\end{array}\right.$$

On $S_1,$ our new $M$-system $[M_{p,J}(x)]=[L_{p}(2x)[M_J(x)]^t]^t$ is the $M$-system $M_{p,J}$ defined by :

$$[M_{p,J}](x)=\left\{\begin{array}{rrrrrrr}
{[p_1(x)],}&\mbox{if}\ \;x\in T_{1,1},\\
{[p_0(x),p_1(x)]},&\mbox{if}\ \;x\in \pm [\frac{1}{14},\frac{1}{7})\cup \pm [\frac{3}{14}, \frac {1}{4}) \subseteq T_{1,2}, \\
{[p_0(x),p_1(x)]},&\mbox{if}\ \;x\in \pm [\frac{1}{4}, \frac{2}{7}) \subseteq T_{1,2}, \\
{[0,\sqrt{2},0],}&\mbox{if}\ \;x\in \pm [\frac{3}{7},\frac{1}{2})= T_{1,3}\backslash T_{2,3},\\
{[p_0(x),0,p_1(x)],}&\mbox{if}\;x\in [-\frac{1}{14},\frac{1}{14})=T_{1,3}\cap T_{2,3},\\
{[p_0(x+\frac{1}{2}),p_1(x+\frac{1}{2})],}&\mbox{if}\;x\in T_{2,2},\\
{[p_0(x+\frac{1}{2}),0,p_1(x+\frac{1}{2})],}&\mbox{if}\;x\in T_{2,3}.
\end{array}\right.$$
  
We thus obtain the generalized filter functions coming from $M_{p,J}:$
$$h^{p}_{1,1}(x)\;=\;p_0(x)\chi_{[-\frac{2}{7},\frac{2}{7})}(x),$$
$$h^p_{2,1}(x)\;=\;\sqrt{2}\chi_{[-\frac{1}{2},-\frac{3}{7})\cup [\frac{3}{7},\frac{1}{2})}(x),$$
$$g^p_1(x)\;=\;p_1(x)\chi_{[-\frac{2}{7}, \frac{1}{7}) \cup [\frac{1}{7}, \frac{2}{7})}(x),$$
$$h^p_{1,2}(x)\;=\;p_0(x+\frac{1}{2})\chi_{[-\frac{1}{7},\frac{1}{7})}(x),$$
$$h^p_{2,2}(x)\;=\; 0,$$
$$g^p_2(x)\;=\;p_1(x+\frac{1}{2})\chi_{[-\frac{1}{7},\frac{1}{7})}(x).$$

Note that these are exactly the filter functions obtained in Example \ref{ExJourne}.

\end{example}

One could no doubt adapt the above example to obtain generalized filters similar to those given in  Examples 4.2 and 4.3 of \cite{BCM}.
   
\begin{acknowledgments}
The authors gratefully acknowledge helpful conversations with  
Astrid An Huef and Iain Raeburn.
\end{acknowledgments}

\bibliographystyle{amsalpha}

\end{document}